\documentclass[12pt,a4paper]{article}
\usepackage{amssymb,amsmath,amsthm,latexsym}
\usepackage{fullpage}
\usepackage{graphicx}
\theoremstyle{plain}
\newtheorem{theorem}{Theorem}[section]
\newtheorem{lemma}[theorem]{Lemma}
\newtheorem{proposition}[theorem]{Proposition}

\theoremstyle{definition}

\begin{document}

\title{A historical law of large numbers for the Marcus-Lushnikov process}
\author{St\'ephanie Jacquot\thanks{University of Cambridge, Statistical Laboratory, Centre for Mathematical Sciences
Wilberforce Road, Cambridge, CB3 0WB, UK
e-mail: smj45@cam.ac.uk}}
\date{}
\maketitle
\begin{footnotesize}
\noindent {\bf Abstract}\\

The Marcus-Lushnikov process is a finite stochastic particle system, in which each particle is entirely characterized by its mass. Each pair of particles with masses $x$ and $y$ merges into a single particle at a given rate $K(x,y)$. Under certain assumptions, this process converges to the solution to Smoluchowski equation, as the number of particles increases to infinity. The Marcus-Lushnikov process gives at each time the distribution of masses of the particles present in the system, but does not retain the history of formation of the particles. In this paper, we set up a historical analogue of the Marcus-Lushnikov process (built according the rules of construction of the usual Markov-Lushnikov process) each time giving what we call the historical tree of a particle. The historical tree of a particle present in the Marcus-Lushnikov process at  a given time $t$ encodes information about the times and masses of the coagulation events that have formed that particle. We prove a law of large numbers for the empirical distribution of such historical trees. The limit is a natural measure on trees which is constructed from a solution to Smoluchowski coagulation equation.

\end{footnotesize}

\section{Presentation of the problem}
\subsection{Introduction}

Let $E=(0,\infty)$. Let $K:E\times E\to(0,\infty)$ be a symmetric
continuous function. Let $S$ be the set of finite integer-valued
measures on $(0,\infty)$. $S$ contains elements of the form
\begin{displaymath}
x=\sum_{i=1}^{n}m_{i}\delta_{y_{i}}
\end{displaymath}
for $n\in \mathbb{N}$ where $y_{1},\ldots,y_{n}>0$ are distinct and
for $i\in \{1,\ldots,n\}, m_{i}\in \mathbb{N}$. The Marcus-Lushnikov
process with coagulation kernel $K$ is the continuous time Markov
chain on $S$ with non-zero transition rates given by
\begin{displaymath}
q(x,x^{\prime}) = \left\{ \begin{array}{ll}
m_{i}m_{j}K(y_{i},y_{j}) & \textrm{if $i<j$ }\\
\frac{1}{2}m_{i}(m_{i}-1)K(y_{i},y_{i}) & \textrm{if $i=j$}
\end{array} \right.
\end{displaymath}
for
$x^{\prime}=x+\delta_{y_{i}+y_{j}}-\delta_{y_{i}}-\delta_{y_{j}}.$

Let us give a way of constructing a Marcus-Lushnikov process $(X_
{t})_{t\geq 0}$. Let $[N]=\{1,\ldots,N\}$ and let
$y_{1},\ldots,y_{N}>0$ be the masses (not necessarily distinct) associated to each particle in
$[N]$. Set
\begin{displaymath}
X_{0}=\sum_{i=1}^{N}\delta_{y_{i}}.
\end{displaymath}
For each $i<j$ take an independent random variable $T_{ij}$ such
that $T_{ij}$ is exponential with parameter $K(y_{i},y_{j})$, and
define
\begin{displaymath}
T=\min_{i<j}T_{ij}.
\end{displaymath}
Set $X_{t}=X_{0}$ for $t<T$ and
\[X_{T}=X_{t}-(\delta_{y_{i}}+\delta_{y_{j}}-\delta_{y_{i}+y_{j}})\] if $T=T_{ij}$, then begin
the construction afresh from $X_{T}$. Each pair of clusters $i<j$
coagulates at rate $K(y_{i},y_{j})$.

Now fix $\mu_{0}$ a measure on $E$ and take a sequence
$\left(x_{0}^{N}\right)_{N \geq 0}\in S$ such that
\begin{equation}
\mu_{0}^{N}=N^{-1}x_{0}^{N} \to \mu_{0}
\end{equation}
weakly on $E$ as $N \to \infty$. Let $(X_{t}^{N})_{t\geq 0}$ be
Marcus-Lushnikov with kernel $\frac{K}{N}$ starting from $x_{0}^{N}$. For each $N$, we can write
\[
x_{0}^{N}=\sum_{i=1}^{m(N)}\delta_{y_{i}^{N}}
\]
with $m(N) \leq N$ and $y_{1}^{N},\ldots,y_{N}^{m(N)}>0$ non necessarily distinct.
Set
\[\mu_{t}^{N}=N^{-1}X_{t}^{N}.\]
Without loss of generality, for the rest of the paper, we will take $m(N)=N$. Indeed, after one step in the process above, the number of masses in the system will be less than $N$ and so we will be exactly in the case where $m(N)<N$.

Our aim in this paper is to set up a historical analogue of the
process $\mu_{t}^{N}$ and to prove that it converges to a limit
measure that can be constructed from the strong solution to a
generalised form of Smoluchowski's equation [1],[2](to be made precise
below).

Before defining precisely this new process let us explain why it is
interesting to know about the history of formation of a cluster.

The Marcus-Lushnikov process [5] describes the stochastic Markov
evolution of a finite system of coalescing particles. It gives at each time the
distribution in masses of the particles present in the system but
does not retain any other information that the particles might
contain. In other words, we lose in part the information contained
in the particles that is their history. Why is it interesting to
know about the history ?  For instance, consider a system of $N$
particles with associated masses
$y_{1},\ldots,y_{N}>0.$ Assume that these particles can only be of
three types say either $A,B$ or $C.$ Allow them to coagulate according to the rules of coagulation of the Marcus-Lushnikov process. Then,
the usual Marcus-Lushnikov will give us at each time the masses of
the particles present in the system but will not be able to tell us
for each particle present at this time how many particles of type
$A,B$ or $C$ this particle contains along with the order of
formation. Our historical measure will give us at each time the
particles formed with their respective masses, the time when they
have formed but also the history of the formation from its beginning
that is the time at which each intermediary particles have formed
along with what they contain. We could think of an other application
in industry in the process of making a certain chemical product. We
can assume that in order to make a certain powder we need to put $N$
ingredients in a specific order and at specific times. Then our
historical measure will allow us to follow the formation of the
powder and to detect if ingredients were put in the wrong order at
the wrong time.

We are now going to review the work of [1] and [2] about the
convergence of $\mu_{t}^{N}$ to the strong solution to
Smoluchowski's equation as we will use this tool to prove our main
result (stated in 1.3).
\subsection{Related work}
Take $\varphi:E\to(0,\infty)$ to be a continuous sublinear function.
We suppose that $\varphi \geq 1$. Assume that the coagulation kernel
can be written as follows:
\begin{equation}
K(x,y)=\tilde{K}(x,y)\varphi(x)\varphi(y)
\end{equation} with $\tilde{K}$ bounded  on $E\times E$. For $\mu$
a non-negative Borel measure on $E$ such that,
\[\int_{E\times E}K(x,y)\mu(dx)\mu(dy)<\infty\]
we define $L(\mu)$ as follows:
\[\langle f,L(\mu)\rangle=\frac{1}{2}\int_{E\times E}\big(f(x+y)-f(x)-f(y)\big)K(x,y)\mu(dx)\mu(dy)\] for all $f$ bounded and measurable.
We consider the following measure-valued form of the Smoluchowski
coagulation equation,
\begin{equation}\label{eq1}
\mu_{t}=\mu_{0}+\int_{0}^{t}L(\mu_{s})ds .
\end{equation}
We admit as a strong local solution any map:
\[t\to \mu_{t}:[0,T)\to \mathcal{M}^{+}\] where $T\in (0,\infty]$ and $\mathcal{M}^{+}$ is the set of non-negative Borel measures on $E$, such that
\begin{enumerate}
\item for all $B\subseteq E$ compact,
\[t \to \mu_{t}(B):[0,T)\to (0,\infty)\] is measurable.
\item \[\int_{0}^{t}\int_{E}\varphi^{2}(x)\mu_{s}(dx)<\infty \] for all $t<T$.
\item for all bounded measurable functions $f$,
for all $t<T$,
\begin{align*}
&\int_{E}f(x)\mu_{t}(dx)=\int_{E}f(x)\mu_{0}(dx)\\
&\qquad+\int_{0}^{t}\int_{E\times
E}\big(f(x+y)-f(x)-f(y)\big)K(x,y)\mu_{s}(dx)\mu_{s}(dy)ds.
\end{align*}
\end{enumerate}
In the case $T=\infty$ we call a strong local solution a strong solution. Assume that
\begin{equation}
\langle{\varphi}^{2},\mu_{0}\rangle
=\int_{E}\varphi^{2}(x)\mu_{0}(\mathrm{d}x)<\infty.
\end{equation}
Then [2] tells us that there exists a unique maximal strong solution
to (3) denoted $(\mu_{t})_{t<T}$ for some $T>0.$ Moreover assume
that
\begin{equation}
\int_{E} f(x)\varphi(x)\mu_{0}^{N}(dx) \to \int_{E}
f(x)\varphi(x)\mu_{0}(dx)
\end{equation}
as $N\to\infty$ for all $f$ bounded and continuous on $E$. Then, for
all $t<T$, for all $f$ continuous and bounded on $E$
\begin{displaymath}
\int_{E} f(x)\varphi(x)\mu_{t}^{N}(dx) \to \int_{E}
f(x)\varphi(x)\mu_{t}(dx)
\end{displaymath}
as $N \to \infty$ in probability, that is
\[\mu_{t}^{N}\to \mu_{t}\]
as $N \to \infty$ weakly in probability with weight function
$\varphi$. Indeed, for $\epsilon>0,$ for all $f$ continuous and
bounded on $E$, for all $t<T$, we can find $\nu>0$, such that for all $N>0,$
\begin{equation}
\mathbb{P}\left(\sup_{s\leq t}\left|\int_{E}
f(x)\varphi(x)\mu_{s}^{N}(dx)- \int_{E}
f(x)\varphi(x)\mu_{s}(dx)\right|> \epsilon \right)\leq \exp(-\nu N).
\end{equation}

\subsection{Our main result}

\paragraph{Assumptions}
~~\\
~~\\
The following assumptions hold for the whole paper.
Let $\mu_{0}$ be a measure on $(0,\infty)$ and take a sequence
$\left(\mu_{0}^{N}\right)_{N \geq 0}$ (as defined in subsection 1.1) such that
(1) is satisfied that is
\[\mu_{0}^{N}=\frac{1}{N}\sum_{k=1}^{N}\delta_{y_{k}^{N}}\to \mu_{0} \mbox{ weakly as }N \to
\infty.\] For convenience we will write $y_{k}^{N}=y_{k}$ for all $k
\in [N]=\{1,\ldots,N\}.$ Let $E=(0,\infty)$. Let $K:E\times
E\to(0,\infty)$ be a symmetric continuous function. Let
$\varphi:E\to(0,\infty)$ be a continuous sublinear function. We
suppose that $\varphi \geq 1$. Assume that (2),(4) and (5) are satisfied and that $K$ is bounded on sets of
the form $[0,M]^{2}$ with $M>0.$ Let $(\mu_{t})_{t<T}$ be the strong
solution to (3) with $T>0.$

\paragraph{Basic notations for Trees}
~~\\
~~\\
For a finite set $J\subset \mathbb{N}$, we write
$\mathbb{T}(J)$ for the smallest set containing $J$, such
that $\{i,j\}\in \mathbb{T}(J)$ whenever $i,j\in
\mathbb{T}(J)$. We refer to elements of $\mathbb{T}(J)$ as
trees. They are finite binary trees with leaves labeled by elements
of $J$. Each $i\in \mathbb{T}(J)$ has a set of leaf labels
$\lambda(i)\subseteq J$ determined by
\[\lambda(i)=\{i\}\] for $i\in J$ and
\[\lambda(\{i,j\})=\lambda(i)\cup\lambda(j)\]
for all $i,j\in \mathbb{T}(J)$. For each $i \in
\mathbb{T}(J)$, $\left| \lambda(i) \right|$ will denote the
number of elements of the set $\lambda(i)$. Let
$n:\mathbb{T}(J)\to \mathbb{N}$ be the counting function,
defined as follows:
\[n(i)=1\]for $i\in J$ and
\[n(\{i,j\})=n(i)+n(j)\]
for all $i,j\in \mathbb{T}(J)$. Finally, for $J \subset
\mathbb{N}$ finite, we define
\[\mathbb{T}^{\star}(J)=\mathbb{T}^{\star}J=\{i\in \mathbb{T}(J):
\vert\lambda(i)\vert=n(i)\}.\] This is the set of rooted binary
trees with leaves labeled by distinct elements of $J$.
\paragraph{The coagulation process on Trees}
~~\\
~~\\
According to the rules of construction of the usual Marcus-Lushnikov
process, we set up an analogue of this process, each time
giving the trees present in the system. We start with the set
of labeled particles $[N]=\{1,\ldots,N\}$ with masses
$y_{1},\ldots,y_{N}>0$.  For $i=\{i_{1},i_{2}\}\in
\mathbb{T}^{\star}[N]$, the mass $y_{i}$ of $i$ is given recursively
by
\[y_{i}=y_{i_{1}}+y_{i_{2}}.\] Each tree particle $i \in \mathbb{T}^{\star}[N]$ has a tree type $\tau(i)\in \mathbb{T},$ where $\mathbb{T}=\mathbb{T}(\{1\}),$ determined by $\tau(i)=1$ for $i \in [N],$ and $\tau(\{i,j\})=\{\tau(i),\tau(j)\}$ for $i,j \in \mathbb{T}^{\star}[N].$
We construct the historical analogue of the Marcus-Lushnikov process
on $\mathbb{T}^{\star}[N]$ as follows. Set $I[0]=[N]$. For $i \in
[N],$ set $S_{i}=0$. There exist trees $i[1],\ldots,i[N-1]\in
\mathbb{T}^{\star}[N]$, coagulation times
$S_{i[1]}<\cdots<S_{i[N-1]}<\infty$ (they are exactly the jump times
in the usual Marcus-Lushnikov process), and subsets $I[1],\ldots,
I[N-1] \subset \mathbb{T}^{\star}[N]$, such that for each $n \in
\{1,\ldots,N-1\}$, $(\lambda(i): i \in I[n])$ is a partition of
$[N]$, and $I[n]$ satisfies the recursive relation
\[I[n]=I[n-1]\backslash \{i_{1}[n],i_{2}[n]\}\cup \{i[n]\}\]
where $i[n]=\{i_{1}[n],i_{2}[n]\}.$

Let us give an example with $N=5$. The drawing below represents one
of the configurations we can obtain from the coagulation process on
trees described above.

\setlength{\unitlength}{1mm}
\begin{picture}(50,60)
\put(28,10){\line(0,1){40}} \put(28,10){\vector(1,0){60}}
  \put(28,10){\circle*{1}}
  \put(28,20){\circle*{1}}
  \put(28,30){\circle*{1}}
  \put(28,40){\circle*{1}}
  \put(28,50){\circle*{1}}
  \put(28,50){\line(2,-1){10}}
  \put(28,40){\line(2,1){10}}
  \put(38,45){\circle*{1}}
\put(38,45){\line(0,-1){35}}
  \put(38,45){\line(2,-1){10}}
  \put(58,15){\circle*{1}}
 \put(58,15){\line(0,-1){5}}
  \put(28,10){\line(6,1){30}}
  \put(28,20){\line(6,-1){30}}
 \put(28,30){\line(2,1){20}}
 \put(48,40){\circle*{1}}
\put(48,40){\line(0,-1){30}} \put(58,15){\line(2,1){20}}
\put(48,40){\line(2,-1){30}} \put(78,25){\circle*{1}}
\put(78,25){\line(0,-1){15}} \put(23,9){5} \put(23,19){4}
\put(23,29){3} \put(23,39){2} \put(23,49){1}
\put(38,48){$i[1]=\{1,2\}$} \put(38,5){$S_{i[1]}$}
\put(48,43){$i[2]=\{\{1,2\},3\}$} \put(48,5){$S_{i[2]}$}
\put(55,17){$i[3]=\{4,5\}$} \put(58,5){$S_{i[3]}$}
\put(80,24){$i[4]=\{\{\{1,2\},3\},\{4,5\}\}$} \put(78,5){$S_{i[4]}$}
\end{picture}

Here, $I[1]=\Big\{\{1,2\},3,4,5\Big\},
I[2]=\Big\{\{\{1,2\},3\},4,5\Big\},
I[3]=\Big\{\{\{1,2\},3\},\{4,5\}\Big\}$ and
$I[4]=\Big\{\{\{\{1,2\},3\},\{4,5\}\}\Big\}.$

\paragraph{The historical measures}
~~\\
~~\\
Fix $t<T.$ Our principal object of interest is a process of empirical
particle measures $\tilde{\mu}_{t}^{N}$ on the space of historical
trees $A(0,t)$ which we shall now define. The space $A(0,t)$ is given by
\begin{displaymath}
A(0,t)=\bigcup_{\tau \in \mathbb{T}}A_{\tau}(0,t)
\end{displaymath}
where $A_{1}(0,t)=(0,\infty)$ and for $\tau=\{\tau_{1},\tau_{2}\}\in
\mathbb{T}$,
\[A_{\tau}(0,t)=\big\{(s,\{\xi_{1},\xi_{2}\}):s\in(0,t),\xi_{1}\in A_{\tau_{1}}(0,s),
\xi_{2}\in A_{\tau_{2}}(0,s)\big\}.\]

Let us illustrate these historical spaces through an example. Take $\tau=\{\{1\},1\}\in \mathbb{T}$. The tree below represents a $\xi=\left(s_{2},\{\xi_{1},\xi_{2}\}\right) \in A_{\tau}(0,t)$ with $\xi_{1}\in A_{\{1\}}(0,s_{2}),\xi_{2}\in A_{1}(0,s_{2})$.\\\\

\setlength{\unitlength}{1mm}
\begin{picture}(60,30)
\put(38,10){\line(0,1){20}} \put(38,10){\vector(1,0){45}}
\put(79,5){$\mbox{time}$} \put(73,5){$t$}
\put(73,10){\line(0,1){20}} \put(58,5){$s_{2}$} \put(48,5){$s_{1}$}
\put(47,27){$\xi_{1}$} \put(57,22){$\xi$} \put(38,10){\circle*{1}}
\put(38,20){\circle*{1}} \put(38,30){\circle*{1}}
\put(38,30){\line(2,-1){10}} \put(38,20){\line(2,1){10}}
\put(48,25){\circle*{1}} \put(48,25){\line(2,-1){10}}
\put(38,10){\line(2,1){20}} \put(48,25){\line(0,-1){15}}
\put(58,20){\line(0,-1){10}} \put(58,20){\circle*{1}}
\put(58,20){\line(1,0){15}} \put(24,9){$\xi_{2}=y_{3}$}
\put(33,19){$y_{2}$} \put(33,29){$y_{1}$}
\end{picture}

Here $\xi_{1}=\left(s_{1},\{y_{1},y_{2}\}\right),\xi_{2}=y_{3}$ and
$\xi=\Big(s_{2},\{(s_{1},\{y_{1},y_{2}\}),y_{3}\}\Big)$.

We equip $A(0,t)$ with its Borel $\sigma$-algebra (we explain in
Appendix 5.1 how to equip $A(0,t)$ with a topology). We define on
$A(0,t)$ the mass function $m:A(0,t) \to (0,\infty)$.

For $\xi \in A_{1}(0,t)=(0,\infty)$, we set \[m(\xi)=\xi.\] Recursively for
$\tau=\{\tau_{1},\tau_{2}\} \in \mathbb{T}$, for
$\xi=(s,\{\xi_{1},\xi_{2}\}) \in A_{\tau}(0,t),$ we set
\[m(\xi)=m(\xi_{1})+m(\xi_{2}).\]
The empirical historical measure $\tilde{\mu}_{t}^{N}$ is given by
\begin{displaymath}
\tilde{\mu}_{t}^{N}=\frac{1}{N}\sum_{i\in I(t)}\delta_{\xi_{t}^{i}}
\end{displaymath}
where $I(t) \subset \mathbb{T}^{\star}[N]$  is the set of trees
present in the system at time $t$. $I(t)$ is given by
\[I(t)=I[n] \mbox{ for } S_{i[n]} \leq t < S_{i[n+1]}.\]
For $i\in [N], \xi_{t}^{i}=y_{i}$ and for $i=\{i_{1},i_{2}\}\in
I(t)$, with $S_{i}=s,$ we set
$\xi_{t}^{i}=(s,\{\xi_{s}^{i_{1}},\xi_{s}^{i_{2}}\})$. As we trace
back the past of a particle we obtain a tree. Observe that this
empirical measure $\tilde{\mu}_{t}^{N}$ and our usual Marcus
Lushnikov process $\mu_{t}^{N}$ (defined in subsection 1.1) are related through
the following equality,
\[\mu_{t}^{N}=\tilde{\mu}_{t}^{N}\circ m^{-1}.\]
We are interested in taking the limit of this empirical measure as
$N\to\infty$. We define the limit measure on $A(0,t)$ as follows.
For $\xi \in A_{1}(0,t)$, we set
\[\tilde{\mu}_{t}(d\xi)=\exp\left(-\int_{0}^{t}\int_{E}K(y,y^{\prime})\mu_{r}(dy^{\prime})dr\right)\mu_{0}(d\xi)\]
where $y=m(\xi)$ and $(\mu_{r})_{r<T}$ is the strong deterministic
solution to (3). Recursively for $\tau=\{\tau_{1},\tau_{2}\}\in
\mathbb{T}$, $\xi=(s,\{\xi_{1},\xi_{2}\})\in A_{\tau}(0,t)$ with
$s< t<T$, we define
\[\tilde{\mu}_{t}(d\xi)=\epsilon(\tau)K(m(\xi_{1}),m(\xi_{2}))\tilde{\mu}_{s}(d\xi_{1})\tilde{\mu}_{s}(d\xi_{2})\exp\left(-\int_{s}^{t}\int_{E}K(y,y^{\prime})\mu_{r}(dy^{\prime})\mathrm{d}r\right)ds\] where $y=m(\xi)$, $\epsilon(\tau)=1$ if $\tau_{1}\neq\tau_{2}$ and $\epsilon(\tau)=\frac{1}{2}$ if $\tau_{1}=\tau_{2}$.
\paragraph{Our main result}
~~\\~~\\
Our aim in this paper is to prove the following result.
\begin{theorem}
For all $t<T,$
\begin{equation}
\tilde{\mu}_{t}^{N}\to \tilde{\mu}_{t}
\end{equation}
weakly on $A(0,t)$ in probability.
\end{theorem}
This theorem is proved is Section 4. Before giving in subsection
1.5, an outline of how we are going to prove Theorem 1.1, we need
first to introduce some more material. The next subsection 1.4 is
dedicated to introduce a coupled family of processes built on the same
probability space $\left(\Omega,\mathcal{F},\mathbb{P}\right)$ (to
be specified below), which we shall see are Marcus-Lushnikov and
which will be really convenient to use for most of our intermediary
proofs.

\subsection{A coupled family of Marcus-Lushnikov processes}
We start with the set of particles $[N]$ with associated masses
$y_{1},\ldots,y_{N}>0$ (that is $\mu_{0}^{N}$). For $J \subseteq
[N]$ we set $\mathbb{T}_{+}^{\star}J=\mathbb{T}^{\star}J\backslash
J$.
\paragraph{A probability space and some random variables}
~~\\
~~\\
Let $\Omega=(0,\infty)^{\mathbb{T}_{+}^{\star}[N]}$. A probability
measure $\mathbb{P}$ on $\Omega$ is defined for
$\omega=(\omega_{i}:i \in \mathbb{T}_{+}^{\star}[N])$ by
\begin{displaymath}
\mathbb{P}(d\omega)=\prod_{i \in
\mathbb{T}_{+}^{\star}[N]}\exp(-\omega_{i})d\omega_{i}.
\end{displaymath}

Let $(U_{i}:i\in\mathbb{T}_{+}^{\star}[N])$ be a family of random
variables on $\Omega$ defined as follows. For
$\omega=(\omega_{i}:i\in\mathbb{T}_{+}^{\star}[N])\in \Omega$, set
$U_{i}(\omega)=\omega_{i}$. Then, under $\mathbb{P}$,
$(U_{i}:i\in\mathbb{T}_{+}^{\star}[N])$ is a family of independent
exponential random variables with parameter $1$. Set $\mathcal{F}=
\sigma(U_{i}: i \in \mathbb{T}_{+}^{\star}[N])$. Let us now define
on $\Omega$ a family of random variables from which we will
construct our coupled family of Marcus-Lushnikov processes. For
$i\in [N]$ set \[S_{i}=0.\] For $\{i,j\} \in \mathbb{T}^{\star}[N]$,
define recursively
\begin{displaymath}
S_{\{i,j\}}=\max (S_{i},S_{j})+\frac{N}{K(y_{i},y_{j})}U_{\{i,j\}}.
\end{displaymath}

\paragraph{The coupled family of processes}
~~\\
~~\\
We are going to build a coupled family of processes
$\left(\left(X_{t}^{J}\right)_{t \geq 0}: J \subseteq [N]\right)$ as
follows. Fix $J \subseteq [N]$ and set $n_{J}=|J|$. Set
$\mathcal{F}_{J}=\sigma(U_{i}:i\in \mathbb{T}_{+}^{\star}J)$ and
$\Omega^{J}=(0,\infty)^{\mathbb{T}_{+}^{\star}J}$. Fix
$\omega=\left(\omega_{i}:i \in \mathbb{T}_{+}^{\star}J \right)\in
\Omega^{J}$. We start off with the set of particles $J_{N}^{0}=J$
with respective masses $(y_{i}: i \in J)$. We set
\[X_{0}^{J}=\sum_{i\in J}\delta_{y_{i}}\]
and we consider
\[v^{0}= \min_{\substack{i,j \in J_{N}^{0}\\ $ \scriptsize{with} $ i\neq j}}S_{\{i,j\}}(\omega).\]
Almost surely, we have $v^{0}=S_{\{i_{0},j_{0}\}}(\omega)$ for some unique
$i_{0},j_{0} \in J_{N}^{0}$ with $i_{0}\neq j_{0}$. We then obtain a
new set of particles
$J_{N}^{1}=J_{N}^{0}\backslash\{i_{0},j_{0}\}\cup\{\{i_{0},j_{0}\}\}$.\\
We set
\[X_{t}^{J}(\omega)=X_{0}^{J}\] for $t<v^{0}$ and
\[X_{v^{0}}^{J}(\omega)=\sum_{i \in J_{N}^{1}}\delta_{y_{i}}.\]
Now, starting from $J_{N}^{1}$, we consider
\[v^{1}= \min_{\substack{i,j \in J_{N}^{1}\\ $ \scriptsize{with} $ i\neq
j}}S_{\{i,j\}}(\omega).\] Almost surely, we have
$v^{1}=S_{\{i_{1},j_{1}\}}(\omega)$ for some unique $i_{1},j_{1} \in
J_{N}^{1}$ with $i_{1}\neq j_{1}$. We obtain a new set of particles
$J_{N}^{2}=J_{N}^{1}\backslash\{i_{1},j_{1}\}\cup\{\{i_{1},j_{1}\}\}$.
We set
\[X_{t}^{J}(\omega)=X_{v^{0}}^{J}(\omega)\] for $v^{0} \leq t<v^{1}$ and
\[X_{v^{1}}^{J}(\omega)=\sum_{i
\in J_{N}^{2}}\delta_{y_{i}}.\] We start again as above with
$J_{N}^{2}$ and so on. The process stops when there is only one
particle left in the system.

Therefore, for each $\omega \in \Omega^{J},$ there exist trees
$i[1],\ldots,i[n_{J}-1]\in \mathbb{T}^{\star}J$, coagulation times
$S_{i[1]}<\cdots<S_{i[n_{J}-1]}<\infty$ (they are
some of the $S_{i}(\omega)'$s), and subsets $J_{N}^{1},\ldots,
J_{N}^{n_{J}-1} \subset \mathbb{T}^{\star}J$, such that for each $n
\in \{1,\ldots,n_{J}-1\}$, $\left(\lambda(i): i \in
J_{N}^{n}\right)$ is a partition of $J$, and $J_{N}^{n}$ satisfies
the recursive relation
\[J_{N}^{n}=J_{N}^{n-1}\backslash \{i_{1}[n],i_{2}[n]\}\cup \{i[n]\}\]
where $i[n]=\{i_{1}[n],i_{2}[n]\}.$ We set
\begin{displaymath}
T_{i}^{J}(\omega)=\left\{\begin{array}{ll} S_{i[n]}&
\textrm{if $i \in J_{N}^{n-1}\backslash J_{N}^{n}$}\\ \infty &
\textrm{if $i=i[n_{J}-1]$}
\\ S_{i}(\omega)&
\textrm{if $i \notin \cup_{n=0}^{n_{J}-1}J_{N}^{n}$}\end{array}
\right.
\end{displaymath}
$T_{i}^{J}$ can be interpreted as the time of death of the tree
particle $i$ in $J$, if this particle was alive in the system at
some time. If $J=[N]$, we set $T_{i}^{[N]}=T_{i}$. The empirical
historical measure on trees $\tilde{\mu}_{t}^{N}$ can be rewritten
as follows
\begin{displaymath}
\tilde{\mu}_{t}^{N}=\frac{1}{N}\sum_{i\in \mathbb{T}^{\star}[N]
}1_{\{S_{i}\leq t < T_{i}\}}\delta_{\xi_{t}^{i}}.
\end{displaymath}
\paragraph{The Marcus Lushnikov property for this family}
\begin{theorem}
Let $J \subseteq [N]$. The process $\left(X_{t}^{J}\right)_{t \geq
0}$ is Marcus-Lushnikov with kernel $\frac{K}{N}$ starting from
\[X_{0}^{J}=\sum_{i\in J}\delta_{y_{i}}.\] In particular,
$\left(X_{t}^{[N]}\right)_{t \geq 0}$ and $\left(X_{t}^{N}\right)_{t
\geq 0}$ have same distribution.
\end{theorem}

{\it Proof of Theorem 1.2} : Starting from
\[X_{0}^{J}=\sum_{i\in J}\delta_{y_{i}},\] it is clear that the first jump has the correct distribution
for Marcus-Lushnikov. Let us now look at the
${k+1}^{\scriptsize{th}}$ jump. We condition on
\[A_{k}=\left(X_{0}^{J}=\sum_{i\in J}\delta_{y_{i}},
X_{s_{1}}^{J}=\sum_{i\in
J_{N}^{1}}\delta_{y_{i}},\ldots,X_{s_{k}}^{J}=\sum_{i\in
J_{N}^{k}}\delta_{y_{i}}\right),\] where for $l \in \{1,\ldots,k\}$,
$J_{N}^{l}$ is the set of particles present on $[s_{l},s_{l+1})$ and
$s_{1}<\cdots<s_{k}$ are the jump times. What are the
transition rates to go from step $k$ to step $k+1$? Conditional on $A_{k}$, we consider the set of
particles $J_{N}^{k}$ and look at
\begin{displaymath}
T^{\prime}=\min_{\substack{i,j \in J_{N}^{k} \\$ \scriptsize{with} $
i\neq j}}S_{\{i,j\}}=\min_{\substack{i,j \in J_{N}^{k} \\$
\scriptsize{with} $ i\neq j}}s_{\{i,j\}}+V_{\{i,j\}}
\end{displaymath}
where for each $i,j \in J_{N}^{k}$ with $i \neq j$, $V_{\{i,j\}}$ is
exponential with parameter $\frac{K(y_{i},y_{j})}{N}$ and,
\begin{displaymath}
s_{\{i,j\}} = \left\{ \begin{array}{ll}
0 & \textrm{if $i,j \in J$ }.\\
s_{i} & \textrm{if $i$ was formed at some $s_{i} \in \{s_{1},\ldots,s_{k}\}$ and $j \in J.$}\\
s_{j} & \textrm{if $j$ was formed at some $s_{j} \in \{s_{1},\ldots,s_{k}\}$ and $i \in J.$}\\
\max(s_{i},s_{j}) & \textrm{if $i$ was formed at some $s_{i} \in \{s_{1},\ldots,s_{k}\}$ and} \\
& \textrm{$j$ was formed at some $s_{j} \in
\{s_{1},\ldots,s_{k}\}\backslash\{s_{i}\}$.}
\end{array} \right.
\end{displaymath}
Thus, for each $i,j \in J_{N}^{k}$ with $i \neq j$, $V_{\{i,j\}}$ has started at time $s_{\{i,j\}}$ and has been
running for a duration of $s_{k}-s_{\{i,j\}}$. Nevertheless, by the
memoryless property for exponential random variables, for $h>0$,
\[\mathbb{P}\left(s_{\{i,j\}}+V_{\{i,j\}}>s_{k}+h | s_{\{i,j\}}+V_{\{i,j\}}>s_{k}\right)=\mathbb{P}\left(V_{\{i,j\}}>h\right).\]
Thus, it is equivalent to add $s_{k}-s_{\{i,j\}}$ to $S_{\{i,j\}}$
and consider that $V_{\{i,j\}}$ starts from $s_{k}$. Hence,
conditional on $A_{k}$, we find exactly
the $T$ we considered in the construction of the the
Marcus-Lushnikov process in 1.1. Therefore, the transition rates are
the same. Thus $\left(X_{t}^{J}\right)_{t \geq 0}$ is
Marcus-Lushnikov.
\begin{flushright}
$\square$
\end{flushright}

For all $t\geq 0$, we set $\mu_{t}^{N,J}=N^{-1}X_{t}^{J}$. Note that by
construction, for all $J \subseteq [N], (\mu_{t}^{N,J})_{t\geq 0}$
is measurable with respect to $\mathcal{F}_{J}=\sigma(U_{i}: i \in
\mathbb{T}_{+}^{\star}J)$. We shall see that the coupled family of
Marcus-Lushnikov processes $\left(\left(\mu_{t}^{N,J}\right)_{t \geq
0}: J \subseteq [N]\right)$ will be useful in most of the
intermediary proofs leading to our main result. Let us now give an outline of the intermediary results we need in
order to prove Theorem 1.1.

\subsection{Outline proof of main result}

Fix $t\in[0,T)$. In order to prove Theorem 1.1, we need to prove that for all $f\in
C_{b}\left(A(0,t)\right)$
\begin{equation}\langle
f,\tilde{\mu}_{t}^{N}\rangle=\int_{A(0,t)}f(\xi)\tilde{\mu}_{t}^{N}(d\xi)
\to \langle
f,\tilde{\mu}_{t}\rangle=\int_{A(0,t)}f(\xi)\tilde{\mu}_{t}(d\xi)
\end{equation}
as $N\to \infty$ in probability. To prove (8), we shall see that it
is sufficient to prove that for all $\tau \in \mathbb{T}$, for all
$f\in C_{b}\left(A_{\tau}(0,t)\right)$
\begin{equation}\label{eq5}
\langle
f,\tilde{\mu}_{t}^{N}\rangle=\int_{A_{\tau}(0,t)}f(\xi)\tilde{\mu}_{t}^{N}(d\xi)
\to \langle
f,\tilde{\mu}_{t}\rangle=\int_{A_{\tau}(0,t)}f(\xi)\tilde{\mu}_{t}(d\xi)
\end{equation}
as $N\to \infty$ in probability. Then we will be able to conclude, using a tightness argument (that is explained in subsection 4.1). To show (9), it is
sufficient to prove that
\begin{equation}\mathbb{E}\left(\langle
f,\tilde{\mu}_{t}^{N}\rangle \right)\to \langle
f,\tilde{\mu}_{t}\rangle \end{equation}
and
\begin{equation}
\mathbb{E}\left(\langle f,\tilde{\mu}_{t}^{N}\rangle ^{2}\right)\to \langle f,\tilde{\mu}_{t}\rangle ^{2}
\end{equation}
as $N\to \infty$. Then,
\[\mathbb{E}\left(\langle f,\tilde{\mu}_{t}^{N}-\tilde{\mu}_{t}\rangle ^{2}\right) \to 0 \] and a fortiori
\[\langle f,\tilde{\mu}_{t}^{N}\rangle \to \langle f,\tilde{\mu}_{t}\rangle \]
as $N\to \infty$ in probability.

Thus, we need to compute $\mathbb{E}\left(\langle
f,\tilde{\mu}_{t}^{N} \rangle\right)$. The paper is set as follows.
Section 2 is dedicated to the computation of
$\mathbb{E}\left(\langle f,\tilde{\mu}_{t}^{N} \rangle\right)$. In
Section 3, we will prove (10) and (11) in order to obtain (9). Finally in
Section 4 we will use the results out obtained in the previous
sections to prove Theorem 1.1.

\section{Intermediary computations}
For the whole section we start with the set of particles
$[N]$ with associated masses $y_{1},\ldots,y_{N}>0$ (that is $\mu_{0}^{N}$) and we consider
the probability space $\left(\Omega,\mathcal{F},\mathbb{P}\right)$
that we have defined in subsection 1.4. We fix $0\leq t <T.$ The aim of this
section is to calculate $\mathbb{E}\left(\langle
f,\tilde{\mu}_{t}^{N} \rangle\right)$ for $f\in
C_{b}(A_{\tau}(0,t))$ with $\tau \in \mathbb{T}$.

\subsection{A finite sum of conditional expectations}
Let us fix $\tau \in \mathbb{T}$ with $n$ leaves that is $n(\tau)=n$
and take $f\in C_{b}(A_{\tau}(0,t))$. We are going to express $\mathbb{E}\left(\langle
f,\tilde{\mu}_{t}^{N} \rangle\right)$ as the expectation of a finite sum of conditional expectations. Recall that for $i \in
\mathbb{T}^{\star}[N]$, $\tau(i)$ denotes the type of the labeled
tree $i$. We can write,
\begin{align*}
\langle f,\tilde{\mu}_{t}^{N}\rangle=&\frac{1}{N}\sum_{i \in
\mathbb{T}^{\star}[N]} f\left(\xi_{t}^{i}\right)1_{\{S_{i}\leq
t < T_{i}\}}\\
=&\frac{1}{N}\sum_{\substack{i \in \mathbb{T}^{\star}[N] \\
$\scriptsize{ with }$ \tau(i)=\tau}}
f\left(\xi_{t}^{i}\right)1_{\{S_{i}\leq t < T_{i}\}}
\end{align*}
as $f$ is supported on $A_{\tau}(0,t)$. For $J \subseteq [N],$ and $k \in \mathbb{N}$ define
\[\mathbb{T}^{\star}_{k}J=\{i\in \mathbb{T}^{\star}J: n(i)=k\}.\]
For any $i\in \mathbb{T}^{\star}_{n}[N]$ with type $\tau$, there are $2^{q(\tau)}$ permutations possible of the particles composing $i$ which
will keep the tree $i$ invariant where $q(\tau)$ is the number of
symmetries in the tree $\tau$. This is given recursively
by $q(1)=0$ and for $\tau=\{\tau_{1},\tau_{2}\}\in \mathbb{T}$,
$q(\tau)=q(\tau_{1})+q(\tau_{2})+1_{\{\tau_{1}=\tau_{2}\}}$.
~~\\
~~\\
{\it Example :} Take $i=\{1,\{2,3\}\}$. This tree has for type $\tau=\{1,\{1\}\}$ and $q(\tau)=1$. The permutations leaving the tree $i$ alike are the identity permutation and the one sending $1$ to itself and exchanging $2$ and $3$. Hence $2^{q(\tau)}$ permutations leave this tree invariant.
~~\\
~~\\
Now fix $i_{0}\in \mathbb{T}^{\star}_{n}[n]$ with type $\tau(i_{0})=\tau$. Define
\[S_{n,N}=\Big\{\mbox{ injections }\sigma : [n]\to[N]\Big\}.\] For $\sigma \in S_{n,N},$ set
$\sigma(1,\ldots,n)=\big(\sigma(1),\ldots,\sigma(n)\big).$ Let
\[
[N]^{\star}_{[n]}=\Big\{(i_{1},\ldots,i_{n})\in [N]^{n}: i_{1},\ldots,i_{n} \mbox{ distinct}\Big\}.
\]
Observe that for each $(i_{1},\ldots,i_{n})\in [N]^{\star}_{[n]}$ there is a unique $\sigma\in S_{n,N}$ such that $\sigma(1,\ldots,n)=(i_{1},\ldots,i_{n}).$ Hence,
\[[N]^{\star}_{[n]}=\Big\{\sigma(1,\ldots,n): \sigma\in S_{n,N}\Big\}.\]
For $\sigma \in S_{n,N}$ define $\sigma(i_{0}) \in \mathbb{T}^{\star}_{n}[N]$ to be the tree obtained from $i_{0}$ by replacing each particle $i\in [n]$ in $i_{0}$ by $\sigma(i).$
~~\\
~~\\
{\it Example :} If $N=5$, $i_{0}=\{1,\{2,3\}\}$ and $\sigma \in S_{3,5}$ is such that $\sigma(1)=3, \sigma(2)=5, \sigma(3)=1$, then
\[\sigma(i_{0})=\{3,\{5,1\}\}.\]
Define
\[I_{0}:[N]^{\star}_{[n]}\to \mathbb{T}^{\star}_{n}[N] \mbox{ by }I_{0}\left(\sigma(1,\ldots,n)\right)=\sigma(i_{0}).\]
Thus, we can write
\begin{align*}
\langle
f,\tilde{\mu}_{t}^{N}\rangle &=\frac{2^{-q(\tau)}}{N}\sum_{\sigma\in S_{n,N}}
f\left(\xi_{t}^{\sigma(i_{0})}\right)1_{\{S_{\sigma(i_{0})}\leq
t < T_{\sigma(i_{0})}\}}\\
&=\frac{2^{-q(\tau)}}{N}\sum_{(i_{1},\ldots,i_{n})\in [N]^{\star}_{[n]}}
f\left(\xi_{t}^{I_{0}(i_{1},\ldots,i_{n})}\right)1_{\{S_{I_{0}(i_{1},\ldots,i_{n})}\leq
t < T_{I_{0}(i_{1},\ldots,i_{n})}\}}.
\end{align*}
Hence, setting
$J(\sigma(i_{0}))=[N]\backslash\{\lambda(\sigma(i_{0}))\}$, we can write
\begin{align*}
\mathbb{E}\left(\langle f,\tilde{\mu}_{t}^{N}\rangle\right)
&=\frac{2^{-q(\tau)}}{N}\sum_{\sigma\in S_{n,N}}
\mathbb{E}\bigg(f\left(\xi_{t}^{\sigma(i_{0})}\right)1_{\{S_{\sigma(i_{0})}\leq
t < T_{\sigma(i_{0})}\}}\bigg)\\
&=\mathbb{E}\left(\frac{2^{-q(\tau)}}{N}\sum_{\sigma\in S_{n,N}}
\mathbb{E}\bigg(f\left(\xi_{t}^{\sigma(i_{0})}\right)1_{\{S_{\sigma(i_{0})}\leq
t < T_{\sigma(i_{0})}\}}\vert
\mathcal{F}^{J(\sigma(i_{0}))}\bigg)\right).
\end{align*}

We shall compute
\begin{equation}\label{eq7}
\mathbb{E}\bigg(f
\left(\xi_{t}^{\sigma(i_{0})}\right)1_{\{S_{\sigma(i_{0})}\leq
t < T_{\sigma(i_{0})}\}}\vert
\mathcal{F}^{J(\sigma(i_{0}))}\bigg)
\end{equation}
for all $\sigma\in S_{n,N}$. In order to compute this quantity, we shall find it useful
to work with labeled trees.

\subsection{Working with labeled trees}
We are going to introduce spaces similar to $A(0,t)$, but for
labeled trees.
\paragraph{Vector tree of masses}~~\\~~\\
Each particle $j \in [N]$ has a given mass $y_{j}.$ Take a tree
$i\in \mathbb{T}_{n}^{\star}[N]$. Assume that this tree is formed from the particles $i_{1}<\cdots<i_{n}$. We write
$y=(y_{i_{1}},\ldots,y_{i_{n}})$ for its associated vector of masses. Its vector tree of masses $\tilde{y}$ is
the tree of masses obtained from $i$ by replacing each
particle in this tree by its mass. For instance $i=\{1,\{2,3\}\}$ and $y=(y_{1},y_{2},y_{3})$ give $\tilde{y}=\{y_{1},\{y_{2},y_{3}\}\}.$

\paragraph{Labeled historical spaces of trees}~~\\~~\\
Define $A_{k}^{y_{k}}(0,t)$ for $k \in \{1,\ldots,N\}$, by
\[A_{k}^{y_{k}}(0,t)=\{k\}\times \{y_{k}\}\] and recursively for
$i=\{i_{1},i_{2}\}\in
\mathbb{T}_{n}^{\star}[N]$ define
\[A_{i}^{y}(0,t)=\big\{\left(s,\{\xi_{1},\xi_{2}\}\right) :s\in(0,t),\xi_{1}\in
A_{i_{1}}^{y^{1}}(0,s), \xi_{2}\in
A_{i_{2}}^{y^{2}}(0,s)\big\}\]
where $y, y^{1}$ and $y^{2}$ are the respective vectors of masses of $i, i_{1}$ and $i_{2}.$
For $i\in \mathbb{T}_{n}^{\star}[N]$ (without any associated masses) define
\begin{displaymath}
A_{i}(0,t)=\bigcup_{y \in
(0,\infty)^{n}}A_{i}^{y}(0,t)
\end{displaymath}
and set
\begin{displaymath}
\tilde{A}(0,t)=\bigcup_{N=1}^{\infty}\bigcup_{i \in
\mathbb{T}^{\star}[N]}A_{i}(0,t).
\end{displaymath}
Observe that when we integrate over
$A_{i}^{y}(0,t)$, we only integrate over the
coagulation times as the masses of the particles are fixed whereas
integrating over $A_{i}(0,t)$ means integrating over the
coagulation times along with the masses.
\paragraph{Reduction of the problem}~~\\~~\\
Each particle $j\in[N]$ has a given mass $y_{j}$.
Take $\tau \in \mathbb{T}$ with $n(\tau)=n$ for some $n \in
\mathbb{N}$. Fix $f \in C_{b}(A_{\tau}(0,t))$. In order to compute
$\mathbb{E}\left(\langle f,\tilde{\mu}_{t}^{N} \rangle\right)$,
without loss of generality, by subsection 2.1, we need to compute
\begin{displaymath}
\mathbb{E}\left(f
\left(\xi_{t}^{i}\right)1_{\{S_{i}\leq t <
T_{i}\}}\vert \mathcal{F}^{J}\right)N^{n-1}
\end{displaymath}
for $i \in \mathbb{T}_{n}^{\star}[n]$ with type $\tau(i)=\tau$ and
associated vector of masses $y=(y_{1},\ldots,y_{n})$ and $J=[N]\backslash
[n]$. The map
\[g_{i}: A_{i}(0,t)\to A_{\tau}(0,t) \] on
forgetting labels is $2^{q(\tau)}$ to $1$. Define \[f_{i}=f\circ
g_{i}.\]
We have $f_{i}\in
C_{b}\left(A_{i}(0,t)\right)$. Let $\xi_{t}^{i}\in A_{\tau}(0,t)$. Then,
\[g_{i}^{-1}(\xi_{t}^{i})=\Big\{\zeta_{t}\in A_{i}(0,t) : g_{i}(\zeta_{t})=\xi_{t}^{i}\Big\}.\]
The set $g_{i}^{-1}(\xi_{t}^{i})$ contains $2^{q(\tau)}$ elements.
Hence, we can write
\begin{displaymath}
f \left(\xi_{t}^{i} \right)
=2^{-q(\tau)}\sum_{\zeta_{t}\in g_{i}^{-1}(\xi_{t}^{i})}f\left(g_{i}(\zeta_{t})\right).
\end{displaymath}
So,
\begin{displaymath}
\mathbb{E}\Big(f \left(\xi_{t}^{i} \right)1_{\{S_{i}\leq t < T_{i}\}} \Big)
=2^{-q(\tau)}\mathbb{E}\Big(\sum_{\zeta_{t}\in g_{i}^{-1}(\xi_{t}^{i})}f\left(g_{i}(\zeta_{t})\right)1_{\{S_{i}\leq t < T_{i}\}} \Big).
\end{displaymath}

If we fix $\zeta_{t}^{0} \in g_{i}^{-1}(\xi_{t}^{i})$, then all the other elements of $g_{i}^{-1}(\xi_{t}^{i})$ can be obtained from $\zeta_{t}^{0}$ by permuting masses between the particles that form a symmetry in the tree $i$. Hence $g_{i}^{-1}(\xi_{t}^{i})$ can be written as a set depending only on $\zeta_{t}^{0}$.
~~\\
~~\\
{\it Example :} Take $i=\{1,2\}$. Then $\tau(i)=\tau=\{1\}.$ Let $\xi_{t}^{i}=(s,\{y,y^{\prime}\})\in A_{\tau}(0,t)$. Then,
\[g_{i}^{-1}(\xi_{t}^{i})=\Big\{(s,\{\{1\}\times\{y\},\{2\}\times\{y^{\prime}\}\}),(s,\{\{1\}\times\{y^{\prime}\},\{2\}\times\{y\}\})\Big\}.\]
It is clear in the set $g_{i}^{-1}(\xi_{t}^{i})$ that we can obtain one particle from the other by exchanging masses between particles $1$ and $2$ which are symmetric in this tree.
~~\\
~~\\
Therefore,
\begin{displaymath}
\mathbb{E}\Big(f \left(\xi_{t}^{i} \right)1_{\{S_{i}\leq t < T_{i}\}} \Big)
=\mathbb{E}\Big(f_{i}\left(\zeta_{t}^{0}\right)1_{\{S_{i}\leq t < T_{i}\}} \Big).
\end{displaymath}
Hence,
\begin{equation}
\mathbb{E}\bigg(\mathbb{E}\Big(f \left(\xi_{t}^{i} \right)1_{\{S_{i}\leq t < T_{i}\}}\vert \mathcal{F}^{J}\Big) \bigg)
=\mathbb{E}\bigg(\mathbb{E}\Big(f_{i}\left(\zeta_{t}^{0}\right)1_{\{S_{i}\leq t < T_{i}\}}\vert \mathcal{F}^{J} \Big)\bigg).
\end{equation}
For convenience, according
to the context and the spaces we consider,  $\xi_{t}^{i}$
will stand for either an element of $A_{\tau}(0,t)$ or
$A_{i}(0,t)$ and more particularly since the masses of
the particles are fixed as an element of
$A_{i}^{y}(0,t)$ where $y$ is the vector of masses of the tree $i$.

Hence it is enough to compute
\begin{equation}\mathbb{E}\left(f_{i}
\left(\xi_{t}^{i}\right) 1_{\{S_{i}\leq t <
T_{i}\}}|\mathcal{F}_{J}\right).
\end{equation}
We are going to compute (14) first for some particular $i$ and $f_{i}$. Then, we will use these intermediary results to
solve the general case $i\in
\mathbb{T}^{\star}_{n}[n]$ and $f_{i}\in
C_{b}\left(A_{i}(0,t)\right)$ and by the relation (13) we will obtain
an expression for $f\in C_{b}\left(A_{\tau}(0,t)\right)$.

\subsection{Case $f=1$ and $\tau =1$}
We take $\tau=1 \in \mathbb{T}$ and $f=1 \in
C_{b}\left(A_{\tau}(0,t)\right).$ The corresponding set of labeled
trees with type 1 is $[N].$  Without loss of generality take
$i=1$. The corresponding $f_{1} \in
C_{b}\left(A_{i}(0,t)\right)$ defined in subsection 2.2 is
$f_{1}=1$. Also, since $i \in [N]$, we have by
definition $S_{i}=0.$ Hence, in this case, we want to
compute $\mathbb{P}(T_{1}>t\mid{\mathcal{F}}_{J_{1}})$ where
$J_{1}=\{2,\ldots,N\}$.
\begin{theorem}\label{simple case; theorem}
Let $J_{1}=\{2,\ldots,N\}$ and
${\mathcal{F}}_{J_{1}}=\sigma\left(U_{i}:i\in
\mathbb{T}_{+}^{\star}J_{1}\right)$. Then,
\[\mathbb{P}(T_{1}>t\mid{\mathcal{F}}_{J_{1}})=\exp \left(-\int_{0}^{t}\int_{E}K(y_{1},y)
{\mu}^{N,J_{1}}_{r}(dy)dr \right)\mbox{ a.s}\]where $T_{1}$ is as defined
previously the time at which particle $1$ dies in
$({\mu}^{N}_{r})_{r<t}$ and $y_{1}$ is the mass of particle $1$.
\end{theorem}
{\it Proof of Theorem 2.1} : Since ${\mu}^{N,J_{1}}_{t}$ is
measurable with respect to $\mathcal{F}_{J_{1}}$, we can write
${\mu}^{N,J_{1}}_{t}=F^{N}_{t}(U_{i}:i\in
\mathbb{T}_{+}^{\star}J_{1})$ where
$F^{N}_{t}:(0,\infty)^{\mathbb{T}_{+}^{\star}J_{1}}\to(0,\infty)$ is
a measurable function. Fix $\omega^{J_{1}}=(\omega^{J_{1}}_{i}:i\in
\mathbb{T}_{+}^{\star}J_{1})\in
(0,\infty)^{\mathbb{T}_{+}^{\star}J_{1}}$ and set
${\mu}^{N,J_{1}}_{t}(\omega^{J_{1}})=F^{N}_{t}(\omega^{J_{1}}_{i}:i\in
\mathbb{T}_{+}^{\star}J_{1})$ (this notation will be kept for the
whole paper).

There exist $r_{1}(\omega^{J_{1}})<\cdots<r_{n}(\omega^{J_{1}})$
such that ${\mu}^{N,J_{1}}_{r}(\omega^{J_{1}})$ is constant on the
interval $[r_{k}(\omega^{J_{1}}),r_{k+1}(\omega^{J_{1}}))$ for
$k=0,\ldots,n$ with the convention $r_{0}(\omega^{J_{1}})=0$ and
$r_{n+1}(\omega^{J_{1}})=t$.

Moreover, we can write
\[T_{1}=G^{N}_{t}\Big((U_{\{1,i\}}:i\in \mathbb{T}^{\star}J_{1}),(U_{i}:i\in \mathbb{T}_{+}^{\star}J_{1})\Big)\] where
$G^{N}_{t}:(0,\infty)^{\mathbb{T}^{\star}J_{1}}\times(0,\infty)^{\mathbb{T}_{+}^{\star}J_{1}}\to\mathbb{R}$
is a measurable function. Hence,
\begin{displaymath}
T_{1}(\omega^{J_{1}})=G^{N}_{t}\Big((U_{\{1,i\}}:i\in
\mathbb{T}^{\star}J_{1}),\omega^{J_{1}}\Big)
\end{displaymath} is a random variable. To obtain $\mathbb{P}(T_{1}>t\mid{\mathcal{F}}_{J_{1}}),$ by Fubini,
it is enough to compute $\mathbb{P}(T_{1}(\omega^{J_{1}})>t)$.

Now observe that
\[\{T_{1}(\omega^{J_{1}})>t\}=\{T_{1}(\omega^{J_{1}})>t\}\cap\{T_{1}(\omega^{J_{1}})>r_{n}(\omega^{J_{1}})\}\cap\ldots\cap\{T_{1}(\omega^{J_{1}})>r_{1}(\omega^{J_{1}})\}\]
because

\[\{T_{1}(\omega^{J_{1}})>t\}\subseteq\{T_{1}(\omega^{J_{1}})>r_{n}(\omega^{J_{1}})\}\subseteq\ldots\subseteq\{T_{1}(\omega^{J_{1}})>r_{1}(\omega^{J_{1}})\}.\]

Hence,
\begin{equation}\label{eq9}
\mathbb{P}\left(T_{1}(\omega^{J_{1}})>t\right)=\prod_{k=0}^{n}\mathbb{P}\left(T_{1}(\omega^{J_{1}})>r_{k+1}(\omega^{J_{1}})|T_{1}(\omega^{J_{1}})>r_{k}(\omega^{J_{1}})\right).
\end{equation}
\\
Let us start by computing the quantity
$\mathbb{P}\left(T_{1}(\omega^{J_{1}})>r_{1}(\omega^{J_{1}})\right).$
Then, we will calculate
$\mathbb{P}\left(T_{1}(\omega^{J_{1}})>r_{k+1}(\omega^{J_{1}})|T_{1}(\omega^{J_{1}})>r_{k}(\omega^{J_{1}})\right)$
for $k\in\{1,\ldots,n\}$.\\\\ For $\{i,j\} \in
\mathbb{T}^{\star}[N]$ set
\[V_{\{i,j\}}=\frac{N}{K(y_{i},y_{j})}U_{\{i,j\}}.\]
Observe that
$\{T_{1}(\omega^{J_{1}})>r_{1}(\omega^{J_{1}})\}=\{S_{\{1,j\}}(\omega^{J_{1}})>r_{1}(\omega^{J_{1}}):j=2,\ldots,N\}$
and that for $j\in\{2,\ldots,N\}$, $S_{\{1,j\}}(\omega^{J_{1}})
=V_{\{1,j\}}.$ Hence,
\begin{align*}\mathbb{P}\left(T_{1}(\omega^{J_{1}})>r_{1}(\omega^{J_{1}})\right) =&\mathbb{P}\left(V_{\{1,j\}}>r_{1}(\omega^{J_{1}}):j=2,\ldots,N\right)\\
=&\prod_{j=2}^{N}\exp\left(-\frac{K(y_{1},y)}{N}r_{1}(\omega^{J_{1}})\right)\\
=&\exp\left(-\sum_{j=2}^{N}\frac{K(y_{1},y)}{N}r_{1}(\omega^{J_{1}})\right)
\end{align*}
by independence of the $(V_{\{1,j\}}:j\in J_{1})$. Now, for all
$r\in[0,r_{1}(\omega^{J_{1}}))$,
\begin{displaymath}\mu_{r}^{N,J_{1}}(\omega^{J_{1}})=\frac{1}{N}\sum_{j=2}^{N}\delta_{y_{j}}.
\end{displaymath}
So,
\[\int_{0}^{r_{1}(\omega^{J_{1}})}\int_{E}K(y_{1},y)\mu_{r}^{N,J_{1}}(\omega^{J_{1}})(dy)dr=
\sum_{j=2}^{N}\frac{K(y_{1},y_{j})}{N}r_{1}(\omega^{J_{1}}).\]

Hence ,
\[\mathbb{P}\left(T_{1}(\omega^{J_{1}})>r_{1}(\omega^{J_{1}})\right)=
\exp\left(-\int_{0}^{r_{1}(\omega^{J_{1}})}\int_{E}K(y_{1},y){\mu}^{N,J_{1}}_{r}(\omega^{J_{1}})(dy)dr\right).\]

Now, let us compute
$\mathbb{P}\left(T_{1}(\omega^{J_{1}})>r_{k+1}(\omega^{J_{1}})\mid
T_{1}(\omega^{J_{1}})>r_{k}(\omega^{J_{1}})\right)$ for
$k\in\{1,\ldots,n\}$. For $k\in\{1,\ldots,n\}$, let
$J_{k}(\omega^{J_{1}})$ be the set of particles present in
$\left({\mu}^{N,J_{1}}_{r}(\omega^{J_{1}})\right)$ for
$r\in[r_{k}(\omega^{J_{1}}),r_{k+1}(\omega^{J_{1}}))$ and assume
that $j_{k}$ is the particle formed at time $r_{k}(\omega^{J_{1}})$.
We can write $j_{k}=\{j_{k}(1),j_{k}(2)\}$ with
$j_{k}(1),j_{k}(2)\in J_{k-1}(\omega^{J_{1}}).$ Observe that
$S_{j_{k}}(\omega^{J_{1}})=r_{k}(\omega^{J_{1}})$. Then,
\begin{align*}
&\mathbb{P}\left(T_{1}(\omega^{J_{1}})>r_{k+1}(\omega^{J_{1}})| T_{1}(\omega^{J_{1}})>r_{k}(\omega^{J_{1}})\right)\\
&\qquad=
\mathbb{P}\left(S_{\{1,j\}}(\omega^{J_{1}})>r_{k+1}(\omega^{J_{1}})\mbox{
for all }j\in J_{k}(\omega^{J_{1}})|
T_{1}(\omega^{J_{1}})>r_{k}(\omega^{J_{1}})\right).
\end{align*}
Observe that $J_{k}(\omega^{J_{1}})=\{j_{k}\}\cup
J_{k-1}(\omega^{J_{1}})\backslash\{j_{k}(1),j_{k}(2)\}$. By the
memoryless property for exponential random variables,
\begin{align*}
&\mathbb{P}\left(T_{1}(\omega^{J_{1}})>r_{k+1}(\omega^{J_{1}})\mid T_{1}(\omega^{J_{1}})>r_{k}(\omega^{J_{1}})\right)\\
&=\mathbb{P}\left(S_{\{1,j\}}(\omega^{J_{1}})>r_{k+1}(\omega^{J_{1}})\mbox{
for all }j\in J_{k} (\omega^{J_{1}})\backslash\{j_{k}\},\right.\\
&\qquad\qquad\left.
S_{\{1,j_{k}\}}(\omega^{J_{1}})>r_{k+1}(\omega^{J_{1}}) \mid
S_{\{1,j\}}(\omega^{J_{1}})>r_{k}(\omega^{J_{1}})\mbox{ for all
}j\in J_{k}(\omega^{J_{1}})\backslash\{j_{k}\}\right).
\end{align*}

Now we can write
$J_{k}(\omega^{J_{1}})\backslash\{j_{k}\}=I_{k}(\omega^{J_{1}})\cup
K_{k}(\omega^{J_{1}})$ where $K_{k}(\omega^{J_{1}})\subset
\{j_{1},\ldots,j_{k-1}\}$ and $I_{k}(\omega^{J_{1}})\subset
\{2,\ldots,N\}$. For $j\in I_{k}(\omega^{J_{1}})$, we have
$S_{\{1,j\}}(\omega^{J_{1}})=V_{\{1,j\}}$, and for $j \in
K_{k}(\omega^{J_{1}})$, there exist $r_{k_{j}}\in
\{r_{1}(\omega^{J_{1}}),\ldots,r_{k-1}(\omega^{J_{1}})\}$ such that
$S_{j}(\omega^{J_{1}})=r_{k_{j}}$. So, $S_{\{1,j\}}(\omega^{J_{1}})=
S_{j}(\omega^{J_{1}})+V_{\{1,j\}}=r_{k_{j}}+V_{\{1,j\}}$. Hence, for
each $j\in J_{k}(\omega^{J_{1}})\backslash\{j_{k}\},$ $S_{\{1,j\}}$
only depends on $V_{\{1,j\}}$. Also,
$S_{j_{k}}(\omega^{J_{1}})=r_{k}(\omega^{J_{1}})$ so
$S_{\{1,j_{k}\}}(\omega^{J_{1}})=r_{k}(\omega^{J_{1}})+
V_{\{1,j_{k}\}}$. Thus, using the independence of
$\left(V_{\{1,j\}},j\in
J_{k}(\omega^{J_{1}})\backslash\{j_{k}\}\right)$ and the memoryless property for exponential random variables, we obtain that :

\begin{align*}&\mathbb{P}(T_{1}(\omega^{J_{1}})>r_{k+1}(\omega^{J_{1}})\mid T_{1}(\omega^{J_{1}})>r_{k}(\omega^{J_{1}}))\\
&\qquad=\mathbb{P}\left(V_{\{1,j_{k}\}}>r_{k+1}(\omega^{J_{1}})-r_{k}(\omega^{J_{1}})\right)\\
&\qquad\prod_{j\in I_{k}(\omega^{J_{1}})}\mathbb{P}\left(V_{\{1,j\}}>r_{k+1}(\omega^{J_{1}})\mid V_{\{1,j\}}>r_{k}(\omega^{J_{1}})\right)\\
&\qquad\prod_{j\in K_{k}(\omega^{J_{1}})}\mathbb{P}\left(V_{\{1,j\}}>r_{k+1}(\omega^{J_{1}})-r_{k_{j}}\mid V_{\{1,j\}}>r_{k}(\omega^{J_{1}})-r_{k_{j}}\right)\\
&\qquad=\mathbb{P}\left(V_{\{1,j_{k}\}}>r_{k+1}(\omega^{J_{1}})-r_{k}(\omega^{J_{1}})\right)\\
&\qquad\prod_{j\in I_{k}(\omega^{J_{1}})}\mathbb{P}\left(V_{\{1,j\}}>r_{k+1}(\omega^{J_{1}})-r_{k}(\omega^{J_{1}})\right)\\
&\qquad\prod_{j\in K_{k}(\omega^{J_{1}})}\mathbb{P}\left(V_{\{1,j\}}>r_{k+1}(\omega^{J_{1}})-r_{k}(\omega^{J_{1}})\right).\\
\end{align*}
Hence,
\begin{align*}
&\mathbb{P}(T_{1}(\omega^{J_{1}})>r_{k+1}(\omega^{J_{1}})\mid T_{1}(\omega^{J_{1}})>r_{k}(\omega^{J_{1}}))\\
&\qquad=\prod_{j\in J_{k}(\omega^{J_{1}})}\mathbb{P}\left(V_{\{1,j\}}>r_{k+1}(\omega^{J_{1}})-r_{k}(\omega^{J_{1}})\right)\\
&\qquad=\prod_{j\in
J_{k}(\omega^{J_{1}})}\exp\left(-\frac{K(y_{1},y_{j})}{N}(r_{k+1}(\omega^{J_{1}})-r_{k}(\omega^{J_{1}}))\right).
\end{align*}

Now, for all $r\in[r_{k}(\omega^{J_{1}}),r_{k+1}(\omega^{J_{1}}))$,
\begin{displaymath}\mu^{N,J_{1}}_{r}(\omega^{J_{1}})=\frac{1}{N}
\sum_{j\in J_{k}(\omega^{J_{1}})}\delta_{y_{j}}.
\end{displaymath}

So,
\[\int_{r_{k}(\omega^{J_{1}})}^{r_{k+1}(\omega^{J_{1}})}\int_{E}K(y_{1},y)\mu^{N,J_{1}}_{r}(\omega^{J_{1}})(dy)dr=
\sum_{j\in
J_{k}(\omega^{J_{1}})}\frac{K(y_{1},y_{j})}{N}\left(r_{k+1}(\omega^{J_{1}})-r_{k}(\omega^{J_{1}})\right).\]
Hence,
\begin{align*}
&\mathbb{P}(T_{1}(\omega^{J_{1}})>r_{k+1}(\omega^{J_{1}})\mid T_{1}(\omega^{J_{1}})>r_{k}(\omega^{J_{1}}))\\
&\qquad=\exp\left(-\int_{r_{k}(\omega^{J_{1}})}^{r_{k+1}(\omega^{J_{1}})}\int_{E}K(y_{1},y){\mu}^{N,J_{1}}_{r}(\omega^{J_{1}})(dy)dr\right).
\end{align*}

and using the equality (15) we obtain :
\[\mathbb{P}(T_{1}(\omega^{J_{1}})>t)=
\exp\left(-\int_{0}^{t}\int_{E}K(y_{1},y){\mu}^{N,J_{1}}_{r}(\omega^{J_{1}})(dy)dr\right).\]
Finally, by Fubini, we get :
\[\mathbb{P}(T_{1}>t\mid \mathcal{F}_{J_{1}})=
\exp\left(-\int_{0}^{t}\int_{E}K(y_{1},y){\mu}^{N,J_{1}}_{r}(\mathrm{d}y)\mathrm{d}r\right).\]as
required.

\begin{flushright}
$\square$
\end{flushright}

\subsection{Case $f$ and $\tau$ general}
Recall that each particle $j\in[N]$ has a given mass $y_{j}>0$.
Take $\tau \in \mathbb{T}$ with $n(\tau)=n$, for some $n \in
\mathbb{N}$ and fix $i\in \mathbb{T}_{n}^{\star}[n]$ with
type $\tau(i)=\tau$ and associated vector of masses $y=(y_{1},\ldots,y_{n}).$ Fix $f
\in C_{b}\left(A_{\tau}(0,t)\right)$ and consider $f_{i}\in
C_{b}(A_{i}(0,t))$ constructed from $f$ as in subsection 2.2. Let $J=[N]\backslash[n]$ and set
$\mathcal{F}_{J}=\sigma(U_{j}:j \in \mathbb{T}_{+}^{\star}J)$.  Our
aim in this subsection is to compute
\[\mathbb{E}_{y,\mu_{0}^{N,J}}\left(f_{i}\left(\xi_{t}^{i}\right)1_{\{S_{i}
\leq t < T_{i}\}}\vert
\mathcal{F}_{J}\right)\] where
$T_{i}$, the time of death of the labeled tree
$i$, has been defined in the subsection 1.4. The notation
$"\mathbb{E}_{y,\mu_{0}^{N,J}}"$ means that the expectation is taken conditional on starting with the
set of particles $[N]$ with respective masses $y_{1},\ldots,y_{N}>0$. This condition can be rewritten as follows : $y=(y_{1},\ldots,y_{n})$ and
$\mu_{0}^{N,J}=\frac{1}{N}\sum_{j \in J}\delta_{y_{j}}$.

\subsubsection{Useful notations}

Fix $0 \leq t<T.$

\paragraph{The space $\Delta(\xi)$}
~~\\

For $k\in [N]$, for $\xi \in A_{k}^{y_{k}}(0,t)$, define
\[\Delta(\xi)=(0,t)\times\{\xi\}.\]
Recursively for $i \in \mathbb{T}_{n}^{\star}[n]$ with associated vector of masses $y$, for  $\xi=(s,\{\xi_{1},\xi_{2}\}) \in
A_{i}^{y}(0,t)$, define
\[\Delta(\xi)=\Delta(\xi_{1})\cup \Delta(\xi_{2})\cup \big([s,t)\cup\{\xi\}\big).\]

\paragraph{The projection map $\Delta(\xi)\to (0,t)$}
~~\\

For $k\in [N]$, for $\xi \in A_{k}^{y_{k}}(0,t)$, define the projection map  $\Pi=\Pi_{\xi}:\Delta(\xi)\to (0,t)$ by for $s\in (0,t)$,
\[\Pi_{\xi}\left((s,\{\xi\})\right)=s.\]
Recursively for $i \in \mathbb{T}_{n}^{\star}[n]$ with associated vector of masses $y$, for  $\xi=(s,\{\xi_{1},\xi_{2}\}) \in
A_{i}^{y}(0,t)$, define the projection map  $\Pi=\Pi_{\xi}:\Delta(\xi)\to (0,t)$ by
\begin{displaymath}
\Pi(u)=\Pi_{\xi}(u)=\left\{\begin{array}{ll} \Pi_{\xi_{1}}(u)&
\textrm{if $u \in \Delta(\xi_{1})$}\\ \Pi_{\xi_{2}}(u)& \textrm{if
$u \in \Delta(\xi_{2})$}
\\ r&
\textrm{if $u=(r,\xi)\in [s,t)\cup\{\xi\}$}\end{array} \right.
\end{displaymath}

\paragraph{The projection map from $\Delta(\xi)$ onto the space of trees of masses}
~~\\

For $k\in [N]$, for $\xi \in A_{k}^{y_{k}}(0,t)$, set
$y_{t}(\xi)=y_{k}$ and define
\[y=y(\xi):\Delta(\xi)\to (0,\infty) \mbox{ by } y(\xi)(u)=y_{u}=y_{k}.\]
Recursively for $i \in \mathbb{T}_{n}^{\star}[n]$ with associated vector of masses $y$, for  $\xi=(s,\{\xi_{1},\xi_{2}\}) \in
A_{i}^{y}(0,t)$, set \[y_{t}(\xi)=\{y_{s}(\xi_{1}),y_{s}(\xi_{2})\}\]
and define
\[y=y(\xi):\Delta(\xi)\to \mathbb{T}(0,\infty)\] where the space $\mathbb{T}(0,\infty)$ is given by
\begin{displaymath}
\mathbb{T}(0,\infty)=\bigcup_{\tau \in \mathbb{T}}\mathbb{T}_{\tau}(0,\infty)
\end{displaymath}
where $\mathbb{T}_{1}(0,\infty)=(0,\infty)$ and for $\tau=\{\tau_{1},\tau_{2}\}\in
\mathbb{T}$,
\[\mathbb{T}_{\tau}(0,\infty)=\big\{y=\{y^{1},y^{2}\}:y_{1}\in \mathbb{T}_{\tau_{1}}(0,\infty),
y_{2}\in \mathbb{T}_{\tau_{2}}(0,\infty)\big\},\]
by
\begin{displaymath}
y(\xi)(u)=\left\{\begin{array}{ll} y(\xi_{1})(u)& \textrm{if $u \in
\Delta(\xi_{1})$}\\ y(\xi_{2})(u)& \textrm{if $u \in
\Delta(\xi_{2})$}
\\ y_{u}=\{y_{s}(\xi_{1}),y_{s}(\xi_{2})\}&
\textrm{if $u=(r,\xi)\in [s,t)\cup\{\xi\}$}\end{array} \right.
\end{displaymath}

\paragraph{The notation $K_{\xi}$}
~~\\

For $\xi \in A_{k}(0,t)$ with $k\in [N]$, set
\[K_{\xi}=1.\]

Recursively for $i\in
\mathbb{T}^{\star}[N]$, for $\xi=(s,\{\xi_{1},\xi_{2}\})\in
A_{i}(0,t)$, define
\[K_{\xi}=K(m(\xi_{1}),m(\xi_{2}))K_{\xi_{1}}K_{\xi_{2}}.\]

\subsubsection{An expression for the conditional expectation}

Set $\mathcal{F}_{[n]}=\sigma(U_{j}:j \in
\mathbb{T}_{+}^{\star}[n]).$ The aim is to compute
\begin{align*}
&\mathbb{E}_{y,\mu_{0}^{N,J}}\left(f_{i}\left(\xi_{t}^{i}\right)1_{\{S_{i}
\leq t < T_{i}\}}\vert \mathcal{F}_{J}\vee
\mathcal{F}_{[n]}\right)\\
&\qquad=f_{i}\left(\xi_{t}^{i}\right)1_{\{S_{i}
\leq t < T_{i}^{[n]}\}}
\mathbb{E}_{y,\mu_{0}^{N,J}}\left(1_{\{t < T_{i}\}}\vert
\mathcal{F}_{J}\vee \mathcal{F}_{[n]}\right)
\end{align*}
because
$f_{i}\left(\xi_{t}^{i}\right)1_{\{S_{i}
\leq t < T_{i}^{[n]}\}}$ is measurable with respect to
$\mathcal{F}_{[n]}$ and $\{t < T_{i}\} \subset \{t <
T_{i}^{[n]}\}$ ($T_{i}^{[n]}$ is defined in
1.4). The lemma below will be proved later on.

\begin{lemma}
With the notations above, on the event $\{S_{i} \leq t <
T_{i}^{[n]}\}$,
\[
\mathbb{P}_{y,\mu_{0}^{N,J}}\left(t < T_{i}\vert
\mathcal{F}_{J}\vee \mathcal{F}_{[n]}\right)
=\exp\left(-\int_{\Delta(\xi_{t}^{i})}\int_{E}K(y_{r},y^{\prime})\mu_{\Pi_{\xi_{t}^{i}}(r)}^{N,J}(dy^{\prime})dr
\right) \mbox{ a.s}.
\]

\end{lemma}

Hence, integrating over $A_{i}^{y}(0,t)$ we
obtain
\begin{align*}
\mathbb{E}_{y,\mu_{0}^{N,J}}\left(f_{i}\left(\xi_{t}^{i}\right)1_{\{S_{i}
\leq t < T_{i}\}}\vert \mathcal{F}_{J}\right)=&\int_{A_{i}^{y}(0,t)}f_{i}\left(\xi\right)\mathbb{P}_{y}\left(S_{i}\leq t < T_{i}^{[n]},\xi_{t}^{i}\in d\xi\right)\\
&\exp\left(-\int_{\Delta(\xi)}\int_{E}K(y_{r},y^{\prime})\mu_{\Pi_{\xi}(r)}^{N,J}(dy^{\prime})dr
\right).
\end{align*}
Therefore, we need to compute
$\mathbb{P}_{y}\left(S_{i}\leq t <
T_{i}^{[n]},\xi_{t}^{i}\in d\xi\right)$ for
some $\xi\in A_{i}^{y}(0,t)$.

\begin{lemma}

\begin{displaymath}
\mathbb{P}_{y}\left(S_{i}\leq t <
T_{i}^{[n]},\xi_{t}^{i}\in d\xi\right)
=\frac{K_{\xi}}{N^{n-1}}\exp\left(-\int_{0}^{t}\frac{K_{s}(\xi)}{N}ds\right)\nu(d\xi)
\end{displaymath}
 where $\nu=\nu_{i}^{y}$ is the Lebesgue measure defined in the Appendix $5.2$ and for $\xi\in A_{i}^{y}(0,t)$ and $s \in (0,t)$,
\[K_{s}(\xi)=\frac{1}{2}\sum_{\substack{r,r^{\prime} \in \Pi_{\xi}^{-1}(s) \\ r \neq r^{\prime}}}K(y_{r},y_{r^{\prime}})\]
with $\Pi_{\xi} : \Delta(\xi) \to (0,t)$ as defined in 2.4.1.

\end{lemma}

{\it Proof of Lemma 2.3:} Let us fix $\xi\in
A_{i}^{y}(0,t)$. To $\xi$ we can uniquely
associate $\{(s_{1},j_{1}),\ldots,(s_{n-1},j_{n-1})\}$ where
$s_{1}<\cdots<s_{n-1}$ are the coagulation times and
$j_{1},\ldots,j_{n-1}$ are the labeled trees (subtrees of
$i$) formed at $s_{1},\ldots,s_{n-1}$ respectively.
Denote by $J^{k}(\xi)$ the set of trees (subtrees of
$i$) present on $[s_{k},s_{k+1})$ for
$k\in\{0,\ldots,n-1\}$ with the convention $s_{0}=0$ and $s_{n}=t$.
For each $k\in\{1,\ldots,n-1\}$, we write
$j_{k}=\{j_{k}(1),j_{k}(2)\}$ with $j_{k}(1),j_{k}(2)\in
J^{k-1}(\xi)$. Observe that $J^{k}(\xi)=\{j_{k}\}\cup
J^{k-1}(\xi)\backslash \{j_{k}(1),j_{k}(2)\}$. For convenience, for
each $k\in\{1,\ldots,n-1\}$ we will write
\[K_{j_{k}}=K(y_{j_{k}(1)},y_{j_{k}(2)}).\]

We start with
\[\mu_{0}^{N,[n]}=\frac{1}{N}\sum_{j=1}^{n}\delta_{y_{j}}.\]
Let $\left(Y_{k}\right)_{k=1,\ldots,n-1}$ be the jump chain
associated to the Markov-Lushnikov process
$\left(\mu_{r}^{N,[n]}\right)_{r<t}$. Let $J_{1},\ldots,J_{n-1}$ be
the jump times. Observe that after the ${n-1}^{\scriptsize{th}}$
jump there is only one particle left in the system, and so the
process is in an absorbing state. Thus, the configuration $\xi$ can
be represented through the jump times and the jump states by, for
$k=0,\ldots,n-1$,
\[Y_{k}=\frac{1}{N}\sum_{j \in J^{k}(\xi)}\delta_{y_{j}}\mbox{
and } J_{k}=s_{k}.\]
For $k \in \{0,\ldots,n-2\}$, the rate of going from
$Y_{k}=\frac{1}{N}\sum_{j \in J^{k}(\xi)}\delta_{y_{j}} \mbox{ to }
Y_{k+1}=\frac{1}{N}\sum_{j \in J^{k+1}(\xi)}\delta_{y_{j}}$ is
$q_{k,k+1}=\frac{K_{j_{k+1}}}{N}$ and the rate of leaving the state
$Y_{k}=\frac{1}{N}\sum_{j \in J^{k}(\xi)}\delta_{y_{j}}$ is
$q_{k}=\frac{1}{2}\sum_{\substack{j,l \in J^{k}(\xi) \\$
\scriptsize{with} $ j\neq l}}\frac{K(y_{j},y_{l})}{N}$. Hence,
\begin{align*}
&\mathbb{P}_{y}\left(S_{i}\leq t <
T_{i}^{[n]},\xi_{t}^{i}\in d\xi\right)\\
\qquad&=\mathbb{P}_{y}\left(J_{1}\in ds_{1},\ldots,J_{n-1}\in
ds_{n-1},Y_{1}=\frac{1}{N}\sum_{j \in
J^{0}(\xi)}\delta_{y_{j}},\ldots,Y_{n-1}=\frac{1}{N}\sum_{j \in
J^{n-1}(\xi)}\delta_{y_{j}}\right)\\ \qquad&=q_{0,1}\ldots
q_{n-2,n-1}\exp\left(-q_{1}s_{1}\right)\ldots
\exp\left(-q_{n-2}(s_{n-1}-s_{n-2})\right)ds_{1}\ldots ds_{n-1}.
\end{align*}
Replacing $q_{k}$ and $q_{k,k+1}$ by their respective values we
obtain that
\begin{align*}
\mathbb{P}_{y}\left(S_{i}\leq t <
T_{i}^{[n]},\xi_{t}^{i}\in d\xi\right)
&=\bigg(\prod_{k=1}^{n-1}\frac{K_{j_{k}}}{N}\bigg)
\exp\left(-\frac{1}{2}\sum_{\substack{j,l \in J^{1}(\xi)
\\ $\scriptsize{ with }$ j\neq
l}}\frac{K(y_{j},y_{l})}{N}s_{1}\right)\times\ldots\\
&\qquad \ldots\times \exp\left(-\frac{1}{2}\sum_{\substack{j,l \in
J^{n-1}(\xi)
\\ $\scriptsize{ with }$ j\neq
l}}\frac{K(y_{j},y_{l})}{N}(s_{n-1}-s_{n-2})\right)\nu_{i}^{y}(d\xi)\\
&=\frac{K_{\xi}}{N^{n-1}}\exp\left(-\int_{0}^{t}\frac{K_{s}(\xi)}{N}ds\right)\nu_{i}^{y}(d\xi)
\end{align*}
as required.
\begin{flushright}
$\square$
\end{flushright}

Now let us prove Lemma 2.2. Before proving it, let us introduce
some basic notations.
For $\omega=\left(\omega_{j}:j \in \mathbb{T}_{+}^{\star}[N]\right)\in \Omega$, we will write $\omega^{[n]}=\left(\omega_{j}:j \in \mathbb{T}_{+}^{\star}[n]\right)$ and $\omega^{J}=\left(\omega_{j}:j \in \mathbb{T}_{+}^{\star}J\right)$.\\

{\it Proof of Lemma 2.2:} Fix $\xi\in
A_{i}^{y}(0,t)$. We can find
$\omega_{0}^{[n]}=(\omega_{j}^{0}: j \in \mathbb{T}_{+}^{\star}[n])$
representing the configuration $\xi$. To $\xi$ we can uniquely
associate $\{(s_{1},j_{1}),\ldots,(s_{n-1},j_{n-1})\}$ where
$s_{1}<\cdots<s_{n-1}$ are the coagulation times and
$j_{1},\ldots,j_{n-1}$ are the labeled trees (subtrees of
$i$) formed at $s_{1},\ldots,s_{n-1}$ respectively.
Denote by $J^{k}(\xi)$ the set of particles from $\xi$ present on
$[s_{k},s_{k+1})$ for $k\in\{0,\ldots,n-1\}$ with the convention
$s_{0}=0$ and $s_{n}=t$. For each $k\in\{1,\ldots,n-1\}$, we write
$j_{k}=\{j_{k}(1),j_{k}(2)\}$ with $j_{k}(1),j_{k}(2)\in
J^{k-1}(\xi)$. Observe that $J^{k}(\xi)=\{j_{k}\}\cup
J^{k-1}(\xi)\backslash \{j_{k}(1),j_{k}(2)\}$. Also, for
convenience, for each $k\in\{1,\ldots,n-1\}$  we will write
\[K_{j_{k}}=K(y_{j_{k}(1)},y_{j_{k}(2)}).\]
Now, fix $\omega_{0}^{J}=(\omega_{i}^{0}:i\in
\mathbb{T}_{+}^{\star}J)$. There exist
$r_{1}(\omega_{0}^{J})<\cdots<r_{n^{\prime}-1}(\omega_{0}^{J})$ such
that $\mu_{r}^{N,J}(\omega_{0}^{J})$ constant on
$[r_{k}(\omega_{0}^{J}),r_{k+1}(\omega_{0}^{J}))$ for
$k=0,\ldots,n^{\prime}-1$ with the convention
$r_{0}(\omega_{0}^{J})=0$ and $r_{n^{\prime}}(\omega_{0}^{J})=t$.
For each $k\in\{0,\ldots,n^{\prime}-1\}$ denote by
$J_{k}(\omega_{0}^{J})$ the set of particles from
$\mu_{r}^{N,J}(\omega_{0}^{J})$ present on
$[r_{k}(\omega_{0}^{J}),r_{k+1}(\omega_{0}^{J}))$. Define
\[\Omega(\omega_{0}^{J},\xi)=\{\omega \in \Omega
:\omega^{[n]}=\omega_{0}^{[n]}, \omega^{J}=\omega_{0}^{J} \}.\] To
obtain our result, conditional on $\{S_{i} \leq t <
T_{i}^{[n]}\}$ it is enough to compute :
\begin{equation}
\mathbb{P}_{y,\mu_{0}^{N,J}}\left(\{\omega \in
\Omega(\omega_{0}^{J},\xi) : t < T_{i}(\omega)\}\right)
=\mathbb{P}_{y,\mu_{0}^{N,J}}\left(t <
T_{i}(\omega_{0}^{J},\omega_{0}^{[n]})\right)
\end{equation}
for short where $T_{i}(\omega_{0}^{J},\omega_{0}^{[n]})$
is understood as a random variable on $\Omega(\omega_{0}^{J},\xi).$
Since we are conditioning on $\{S_{i} \leq t <
T_{i}^{[n]}\}$, on each interval $[s_{k},s_{k+1}),$ the
particles present at time $s_{k}$ in $\xi$ that is $J^{k}(\xi)$,
will not interfere between each other until time $s_{k+1}$. Hence,
conditional on $\{S_{i} \leq t <
T_{i}^{[n]}\}$, the quantity (16) means that on each
interval $[s_{k},s_{k+1}),$ each particle in $J^{k}(\xi)$ will not
interfere with the configuration
$\left(\mu_{r}^{N,J}(\omega_{0}^{J})\right)_{r\geq s_{k}}$ until
$s_{k+1}$. This is equivalent to consider on $[s_{k},s_{k+1}),$  for
each $i_{0}\in J^{k}(\xi),$ starting from
\[\mu_{0}^{J_{k}(s_{k})\cup \{i_{0}\}}=\sum_{j \in J_{k}(s_{k})\cup \{i_{0}\}}\delta_{y_{j}},\] that $i_{0}$ will survive until $s_{k+1}$ where
$J_{k}(s_{k})$ is the set of particles present in
$\mu_{s_{k}}^{N,J}(\omega_{0}^{J})$ (this set is one of the $J_{p}(\omega_{0}^{J})$ for $p\in\{0,\ldots,n^{\prime}-1\}$). If $T_{i_{0}}^{J_{k}(s_{k})\cup
\{i_{0}\}}$ denotes the time of death of the particle $i_{0}$ when starting
from the set of particles $J_{k}(s_{k})\cup \{i_{0}\}$ (this death time is defined in 1.4) we want $T_{i_{0}}^{J_{k}(s_{k})\cup \{i_{0}\}}(\omega)>s_{k+1}$ for
all $\omega \in \Omega(\omega_{0}^{J},\xi)$. It is
exactly the result we have obtained in Theorem 2.1. For short
we will write, $T_{i_{0}}^{J_{k}(s_{k})\cup
\{i_{0}\}}=T_{i_{0}}^{s_{k}}.$ Thus, setting
\begin{align*}
A_{1}(\omega_{0}^{J},\xi)&=\{\omega \in \Omega(\omega_{0}^{J},\xi) :
\forall i\in
J^{0}(\xi),T_{i}^{s_{1}}(\omega)>s_{1}\}\\ &\qquad\vdots\\
A_{n-1}(\omega_{0}^{J},\xi)&=\{\omega \in \Omega(\omega_{0}^{J},\xi)
: \forall i\in
J^{n-2}(\xi),T_{i}^{s_{n-1}}(\omega)>s_{n-1}\}\\
A_{n}(\omega_{0}^{J},\xi)&=\{\omega \in \Omega(\omega_{0}^{J},\xi):
T_{j_{n-1}}^{t}(\omega)>t\},
\end{align*}
we can write, conditional on $\{S_{i} \leq t <
T_{i}^{[n]}\}$,
\[
\mathbb{P}_{y,\mu_{0}^{N,J}}\left(t <
T_{i}(\omega_{0}^{J},\omega^{[n]})\right)
=\mathbb{P}_{y,\mu_{0}^{N,J}}\left(A_{1}(\omega_{0}^{J},\xi)\cap\ldots\cap
A_{n}(\omega_{0}^{J},\xi)\right).
\]

The process we consider is a continuous-time Markov process whose
state space is the space of measures on the set $(0,\infty)$. Hence,
fixing the coagulation times $s_{1}<\cdots<s_{n-1}$ we can consider
the process independently on each interval $[s_{k},s_{k+1})$
starting respectively from
$(y^{k},\mu_{s_{k}}^{N,J}(\omega_{0}^{J}))$ where
$y^{0}=(y_{1},\ldots,y_{n})$ and for $k\in\{0,\ldots,n-1\}$,
$y^{k}=(y_{j},j \in J^{k}(\xi)).$ Thus,
\begin{align*}
\mathbb{P}_{y,\mu_{0}^{N,J}}\left( t <
T_{i}(\omega_{0}^{J},\omega_{0}^{[n]})\right)
&=\mathbb{P}_{y,\mu_{0}^{N,J}(\omega_{0}^{J})}\left(A_{1}(\omega_{0}^{J},\xi)\right)\times
\mathbb{P}_{y^{1},\mu_{s_{1}}^{N,J}(\omega_{0}^{J})}\left(A_{2}(\omega_{0}^{J},\xi)\right)\times\ldots
\\
&\qquad\qquad\ldots\times
\mathbb{P}_{y^{n-1},\mu_{s_{n-1}}^{N,J}(\omega_{0}^{J})}\left(A_{n}(\omega_{0}^{J},\xi)\right).
\end{align*}

First, let us compute
\[\mathbb{P}_{y,\mu_{0}^{N,J}(\omega_{0}^{J})}\left(A_{1}(\omega_{0}^{J},\xi)\right)=
\mathbb{P}\left(\{\omega \in \Omega(\omega_{0}^{J},\xi) : \forall
i\in J^{0}(\xi),T_{i}^{J}(\omega)>s_{1}\}\right).\]

Take $i_{0}\in J^{0}(\xi)$. For $j\in J_{k}(\omega_{0}^{J})$ for $k
\in \{1,\ldots,n^{\prime}\}$ either
$S_{\{i_{0},j\}}=V_{\{i_{0},j\}}$ if $j$ is a simple particle, or
there exist $p\in \{0,\ldots,k\}$ such that
$S_{j}(\omega_{0}^{J})=r_{p}(\omega_{0}^{J})$ if $j$ is a composed
particle. In this case
$S_{\{i_{0},j\}}=r_{p}(\omega_{0}^{J})+V_{\{i_{0},j\}}$. Hence, when
$i_{0}\in \{1,\ldots,n\}$,
$T_{i_{0}}^{J}(\omega_{0}^{J},\omega_{0}^{[n]})$ only depends on
$(V_{\{i_{0},j\}}:j\in
J_{k}(\omega_{0}^{J}),k=0,\ldots,n^{\prime}-1)$. Hence,
$\left(T_{j}^{J}(\omega_{0}^{J},\omega_{0}^{[n]}):j\in
J^{0}(\xi)\right)$ are independent random variables. Thus,
\[
\mathbb{P}_{y,\mu_{0}^{N,J}(\omega_{0}^{J})}\left(A_{1}(\omega_{0}^{J},\xi)\right)=\prod_{j\in
J^{0}(\xi)}
\mathbb{P}_{y,\mu_{0}^{N,J}(\omega_{0}^{J})}\left(\{\omega \in
\Omega(\omega_{0}^{J},\xi) : T_{j}^{J}(\omega)>s_{1}\}\right).
\]
For $j\in J^{0}(\xi)$, by Theorem 2.1,
\[
\mathbb{P}_{y,\mu_{0}^{N,J}(\omega_{0}^{J})}(\{\omega \in
\Omega(\omega_{0}^{J},\xi) : T_{j}^{J}(\omega)>s_{1}\})
=\exp\left(-\int_{0}^{s_{1}}\int_{E}K(y_{j},y^{\prime})\mu_{r}^{N,J}(\omega_{0}^{J})
(dy^{\prime})dr\right).
\]

Hence,
\[
\mathbb{P}_{y,\mu_{0}^{N,J}(\omega_{0}^{J})}\left(A_{1}(\omega_{0}^{J},\xi)\right)=
\exp \left(-\int_{0}^{s_{1}}\int_{E}\sum_{j\in
J^{0}(\xi)}K(y_{j},y^{\prime})\mu_{r}^{N,J}(\omega_{0}^{J})(dy^{\prime})dr\right).
\]

Now let us compute
$\mathbb{P}_{y^{k},\mu_{s_{k}}^{N,J}(\omega_{0}^{J})}\left(A_{k}(\omega_{0}^{J},\xi)\right)$
for $k \in \{2,\ldots,n-1\}$. We are starting at time $s_{k-1}$ with
the set of particles $J^{k-1}(\xi)$. By the memoryless property for
exponentially distributed random variables, it is equivalent to
start off at time $0$ with the set $J^{k-1}(\xi)$ and to go to
$s_{k}-s_{k-1}$. Hence, the computation of
$\mathbb{P}_{y^{k},\mu_{s_{k}}^{N,J}(\omega_{0}^{J})}(A_{k}(\omega_{0}^{J},\xi))$
is the same than for
$\mathbb{P}_{y,\mu_{0}^{N,J}(\omega_{0}^{J})}\left(A_{1}(\omega_{0}^{J},\xi)\right)$
replacing $J^{0}(\xi)$ by $J^{k-1}(\xi)$, $j_{1}$ by $j_{k}$, and
$s_{1}$ by $s_{k}-s_{k-1}$. Thus,

\[
\mathbb{P}_{y^{k},\mu_{0
s_{k}}^{N,J}(\omega_{0}^{J})}(A_{k}(\omega_{0}^{J},\xi))= \exp
\left(-\int_{s_{k-1}}^{s_{k}}\int_{E}\sum_{j\in
J^{k-1}(\xi)}K(y_{j},y^{\prime})\mu_{r}^{N,J}(\omega_{0}^{J})(dy^{\prime})dr\right).
\]
Thus,
\begin{align*}
&\mathbb{P}_{y,\mu_{0}^{N,J}}\left(\{\omega \in
\Omega(\omega_{0}^{J},\xi) : t < T_{i}(\omega)\}\right)\\
&\qquad= \prod_{k=1}^{n} \exp
\left(-\int_{s_{k-1}}^{s_{k}}\int_{E}\sum_{j\in
J^{k-1}(\xi)}K(y_{j},y^{\prime})\mu_{r}^{N,J}(\omega_{0}^{J})(dy^{\prime})dr\right)\\
&\qquad=\exp
\left(\sum_{k=1}^{n}-\int_{s_{k-1}}^{s_{k}}\int_{E}\sum_{j\in
J^{k-1}(\xi)}K(y_{j},y^{\prime})\mu_{r}^{N,J}(\omega_{0}^{J})(dy^{\prime})dr\right)\\
&\qquad=\exp\left(-\int_{\Delta(\xi)}\int_{E}K(y_{r},y^{\prime})\mu_{\Pi(r)}^{N,J}(dy^{\prime})dr
\right).
\end{align*}
Hence,
\[\mathbb{P}_{y,\mu_{0}^{N,J}}\left(t < T_{i}\vert \mathcal{F}_{J}\vee \mathcal{F}_{[n]}\right)
=\exp\left(-\int_{\Delta(\xi_{t}^{i})}\int_{E}K(y_{r},y^{\prime})\mu_{\Pi_{\xi_{t}^{i}}(r)}^{N,J}(dy^{\prime})dr
\right) \mbox{ a.s}\] as required.
\begin{flushright}
$\square$
\end{flushright}

Hence,
\begin{align*}
&\mathbb{E}_{y,\mu_{0}^{N,J}}\left(f_{i}\left(\xi_{t}^{i}\right)1_{\{S_{i}\leq
t<T_{i}\}}\mid \mathcal{F}_{J}\right)N^{n-1}\\
&=\int_{A_{i}^{y}(0,t)}f_{i}\left(\xi\right)
K_{\xi}
\exp\left(-\int_{\Delta(\xi)}\int_{E}K(y_{r},y^{\prime})\mu_{\Pi_{\xi}(r)}^{N,J}(dy^{\prime})dr
\right) \exp\left(-\int_{0}^{t}\frac{K_{s}(\xi)}{N}ds
\right)\nu_{i}^{y}(d\xi)\\
&=\int_{A_{i}^{y}(0,t)}f_{i}\left(\xi\right)K_{\xi}
h_{N}\left(\left(\mu_{r}^{N,J}\right)_{r<t},\xi \right)\nu_{i}^{y}(d\xi)
\end{align*}
where \[h_{N}\left(\left(\mu_{r}^{N,J}\right)_{r<t},\xi
\right)=\exp\left(-\int_{\Delta(\xi)}\int_{E}K(y_{r},y^{\prime})\mu_{\Pi_{\xi}(r)}^{N,J}(dy^{\prime})dr
\right) \exp\left(-\int_{0}^{t}\frac{K_{s}(\xi)}{N}ds \right)\] and
where for $y\in (0,\infty)^{n},$ the Lebesgue measure $\nu_{i}^{y}$ is defined in the Appendix 5.2.

\subsubsection{Convergence for teh conditional expectation}
Fix $y=(y_{1},\ldots,y_{n})\in(0,\infty)^{n}$. This is the vector of masses associated to the particles $\{1,\ldots,n\}.$ Let $\tilde{y}$ be the tree vector of masses associated to $i\in \mathbb{T}_{n}^{\star}[n].$
Define
\begin{align*}
P_{t}f_{i}(y)=P_{t}f_{i}(\tilde{y})=&\int_{A_{i}^{y}(0,t)}f_{i}\left(\xi
\right)K_{\xi}\exp\left(-\int_{\Delta(\xi)}\int_{E}K(y_{r},y^{\prime})\mu_{\Pi_{\xi}(r)}
(dy^{\prime})dr\right)\nu_{i}^{y}(d\xi)\\
=&\int_{A_{i}^{y}(0,t)}f_{i}\left(\xi
\right)K_{\xi}h\left(\left(\mu_{r}\right)_{r<t},\xi\right)\nu_{i}^{y}(d\xi)
\end{align*}
where
\[h\left(\left(\mu_{r}\right)_{r<t},\xi\right)=\exp\left(-\int_{\Delta(\xi)}\int_{E}K(y_{r},y^{\prime})\mu_{\Pi_{\xi}(r)}
(dy^{\prime})dr\right).\]

Set
\[P_{t}^{N}f_{i}(y)=P_{t}^{N}f_{i}(\tilde{y})=\mathbb{E}_{y,\mu_{0}^{N,J}}\left(f_{i}\left(\xi_{t}^{i}\right)1_{\{S_{i}\leq t<T_{i}\}}\mid \mathcal{F}_{J}\right)N^{n-1}.\]

\begin{theorem}
\[P_{t}^{N}f_{i}(y)\to P_{t}f_{i}(y) \mbox{ in probability as } N\to\infty .\]
\end{theorem}

{\it Proof of Theorem 2.4} : For a given $\xi \in
A_{i}^{y}(0,t)$, we can write
\begin{align*}
h_{N}\left(\left(\mu_{r}^{N,J}\right)_{r<t},\xi \right)=
&\prod_{k=0}^{n-1}\exp\left(-\frac{1}{2}\sum_{\substack{j,l \in
J^{k}(\xi) \\ $ \scriptsize{with} $ j \neq l}}
\frac{K(y_{j},y_{l})}{N}(s_{k+1}-s_{k})\right)\\
&\prod_{k=0}^{n-1}\exp\left(-\int_{s_{k}}^{s_{k+1}}\int_{E}
\sum_{j\in
J^{k}(\xi)}K(y_{j},y^{\prime})\mu_{r}^{N,J}(dy^{\prime})dr \right)
\end{align*}
where $s_{1}<\cdots<s_{n-1}$ are the coagulation times associated to
$\xi$ and $j_{1},\ldots,j_{n-1}$ are its labeled subtrees (subtrees of $i$) formed at
these respective times. $J^{k}(\xi)$ represents the set of tree particles
(subtrees of $i$) from $\xi$ present on
$[s_{k},s_{k+1})$ for $k=0,\ldots,n-1$ with the convention $s_{0}=0$
and $s_{n}=t.$ As $N\to \infty$, for each $k\in\{1,\ldots,n-1\},$
\[\prod_{k=0}^{n-1}\exp\left(-\frac{1}{2}\sum_{\substack{j,l \in J^{k}(\xi) \\ $ \scriptsize{with} $ j\neq l}} \frac{K(y_{j},y_{l})}{N}(s_{k+1}-s_{k})\right)\to 1.\] Now we want to show that

\begin{align*}
&\prod_{k=0}^{n-1}\exp\left(-\int_{s_{k}}^{s_{k+1}}\int_{E}\sum_{j\in J^{k}(\xi)}K(y_{j},y^{\prime})\mu_{r}^{N,J}(dy^{\prime})dr\right)\\
&\qquad\to
\prod_{k=0}^{n-1}\exp\left(-\int_{s_{k}}^{s_{k+1}}\int_{E}\sum_{j\in
J^{k}(\xi)}K(y_{j},y^{\prime})\mu_{r}(dy^{\prime})dr\right)
\end{align*}
as $N\to \infty$ in probability. It is sufficient to show that for each
$k\in \{1,\ldots,n-1\}$ and $j \in J^{k}(\xi)$,
\[\exp\left(-\int_{s_{k}}^{s_{k+1}}\int_{E}K(y_{j},y^{\prime})\mu_{r}^{N,J}(dy^{\prime})dr\right)\to
\exp\left(-\int_{s_{k}}^{s_{k+1}}\int_{E}K(y_{j},y^{\prime})\mu_{r}(dy^{\prime})dr\right)\]
as $N\to \infty$ in probability. For $k\in\{1,\ldots,n-1\}$ and
$j\in J^{k}(\xi)$ consider the random variable
\[X_{N}^{k,j}=\int_{s_{k}}^{s_{k+1}}K(y_{j},y^{\prime})\mu_{r}^{N,J}(dy^{\prime})dr-\int_{s_{k}}^{s_{k+1}}K(y_{j},y^{\prime})\mu_{r}(dy^{\prime})dr.\]
On $[s_{k},s_{k+1})$,
\begin{displaymath}
\mu_{r}^{N}=\mu_{r}^{N,J}+\frac{1}{N}\sum_{l\in
J_{k}(\xi)}\delta_{y_{l}}.
\end{displaymath}
Hence ,
\begin{align*}
X_{N}^{k,j}&=\int_{s_{k}}^{s_{k+1}}K(y_{j},y^{\prime})\mu_{r}^{N}(dy^{\prime})dr\\
&\qquad-\int_{s_{k}}^{s_{k+1}}K(y_{j},y^{\prime})\mu_{r}(dy^{\prime})dr-
(s_{k+1}-s_{k})\sum_{l\in J_{k}(\xi)}\frac{K(y_{j},y_{l})}{N}.
\end{align*}
Now, we know that
\begin{displaymath}
K(y_{j},y)=\tilde{K}(y_{j},y)\varphi(y_{j})\varphi(y).
\end{displaymath}
For $y\in (0,\infty)$ define $f_{y_{j}}:(0,\infty)\to(0,\infty)$
such that
\begin{displaymath}
f_{y_{j}}(y)=\tilde{K}(y_{j},y)\varphi(y_{j}).
\end{displaymath}
By hypothesis, this function is continuous and bounded. Hence, we
can write
\begin{displaymath}
X_{N}^{k,j}=\int_{s_{k}}^{s_{k+1}}\langle
f_{y_{j}}\varphi,\mu_{r}^{N}\rangle dr-\int_{s_{k}}^{s_{k+1}}
\langle f_{y_{j}}\varphi,\mu_{r}\rangle dr-(s_{k+1}-s_{k})\sum_{l\in
J_{k}(\xi)}\frac{K(y_{j},y_{l})}{N}.
\end{displaymath}
Hence ,
\begin{displaymath}
\left|X_{N}^{k,j}\right|\leq \left|(s_{k+1}-s_{k})\sum_{l\in
J_{k}(\xi)}\frac{K(y_{j},y_{l})}{N}\right| +
(s_{k+1}-s_{k})\sup_{r\in[s_{k},s_{k+1})}\left| \langle
f_{y_{j}}\varphi,\mu_{r}^{N}\rangle - \langle
f_{y_{j}}\varphi,\mu_{r}\rangle \right|.
\end{displaymath}
The first term converges to $0.$ By (6), the second converges to $0$
in probability as $N\to \infty$. Hence $X_{N}^{k,j}\to 0$ in
probability as $N\to \infty$. Hence, for all $\xi \in
A_{i}^{y}(0,t)$,
\begin{displaymath}
h_{N}\Big(\left(\mu_{r}^{N,J}\Big)_{r<t},\xi \right) \to
h\Big(\left(\mu_{r}\right)_{r<t},\xi \Big)
\end{displaymath}
in probability as $N \to \infty$.

Now, on $A_{i}^{y}(0,t)$ $K_{\xi}$ is bounded,
say by some $C>0$. Thus,
\[\left| P_{t}^{N}f_{i}(y)-P_{t}f_{i}(y) \right| \leq
C{\|f_{i}\|}_{\infty}\int_{A_{i}^{\tilde{y}}(0,t)}\left|h_{N}\left(\left(\mu_{r}^{N,J}\right)_{r<
t},\xi \right)-h\Big(\left(\mu_{r}\right)_{r< t},\xi \Big) \right|
d\xi\] Moreover, for all $\xi \in
A_{i}^{y}(0,t)$,
\[\left|h_{N}\left(\left(\mu_{r}^{N,J}\right)_{r<t},\xi \right)\right| \leq 1.\] Hence by the Bounded Convergence Theorem,
\[\mathbb{E}\left(\int_{A_{i}^{y}(0,t)}\left|h_{N}\left(\left(\mu_{r}^{N,J}\right)_{r\leq t},\xi \right)-
h\Big(\left(\mu_{r}\right)_{r\leq t},\xi \Big) \right| d\xi
\right)\to 0\] as $N\to \infty.$ Thus,
\[\mathbb{E}\left(\left| P_{t}^{N}f_{i}(y)-P_{t}f_{i}(y) \right|\right)\to 0\] as $N\to \infty$ and a fortiori
\[P_{t}^{N}f_{i}(y) \to P_{t}f_{i}(y)\] as $N\to \infty$ in probability as required.
\begin{flushright}
$\square$
\end{flushright}

\section{Some convergence results}

The aim of this section is to prove (9, that is, for all $\tau \in \mathbb{T}$, for all $f \in C_{b}\left(A_{\tau}(0,t)\right)$, 

\[
\langle f,\tilde{\mu}_{t}^{N} \rangle\to \langle f,\tilde{\mu}_{t} \rangle
\]
as $N \to \infty$ in probability.

\subsection{Convergence of the expectation}
We recall that we are working with the set of initial particles
$[N]=\{1,\ldots,N\}$ with associated masses $y_{1},\ldots,y_{N}>0$.
Let $\tau \in \mathbb{T}$ with $n$ leaves, that is $n(\tau)=n$. Take
$i \in \mathbb{T}_{n}^{\star}[n]$ with type $\tau$. Associate the map $I:[N]_{[n]}^{\star}\to \mathbb{T}_{n}^{\star}[N]$ as defined in 2.1. For $(i_{1},\ldots,i_{n})\in [N]_{[n]}^{\star}$ with associated masses $y_{i_{1}},\ldots,y_{i_{n}}>0$ we will write $y_{(i_{1},\ldots,i_{n})}$ for the vector of masses and $\tilde{y}_{(i_{1},\ldots,i_{n})}$ for the vector tree of masses associated to $I(i_{1},\ldots,i_{n})$. Let us fix $f \in C_{b}\left(A_{\tau}(0,t)\right)$. We aim to prove
the following result.
\begin{proposition}
\[\mathbb{E}\left(\langle f,\tilde{\mu}_{t}^{N} \rangle \right)\to \langle f,\tilde{\mu}_{t} \rangle\] as $N\to\infty$.
\end{proposition}
{\it Proof of Proposition 3.1 : } For $M>0$, define $\Psi_{M}:
(0,\infty) \to (0,\infty)$ to be

\begin{displaymath}
\Psi_{M}(x)=\left\{\begin{array}{ll} 1& \textrm{if $x \leq M$}\\
-x+(M+1)& \textrm{if $ M \leq x \leq M+1$}
\\ 0&
\textrm{if $x \geq M+1$}\end{array} \right.
\end{displaymath}
We can write,
\[\langle f,\tilde{\mu}_{t}^{N} \rangle=\langle f\Psi_{M}(m),\tilde{\mu}_{t}^{N} \rangle +
\langle f(1-\Psi_{M}(m)),\tilde{\mu}_{t}^{N} \rangle\] where $m :
A(0,t) \to (0,\infty)$ is the mass function defined in subsection 1.3. Hence,
\begin{equation}
\mathbb{E}\left(\langle f,\tilde{\mu}_{t}^{N}
\rangle\right)=\mathbb{E}\left(\langle
f\Psi_{M}(m),\tilde{\mu}_{t}^{N} \rangle\right) +
\mathbb{E}\left(\langle f(1-\Psi_{M}(m)),\tilde{\mu}_{t}^{N}
\rangle\right).
\end{equation}
We use the following lemma that we shall prove below.
\begin{lemma}
\[\mathbb{E}\left(\langle f\Psi_{M}(m),\tilde{\mu}_{t}^{N} \rangle \right)\to \langle f\Psi_{M}(m),\tilde{\mu}_{t} \rangle\] as $N\to\infty$.
\end{lemma}
Thus, taking the limsup and liminf over $N$ in the expression
(17) and applying Lemma 3.2, we obtain,
\begin{equation}\limsup_{N}\mathbb{E}\left(\langle f,\tilde{\mu}_{t}^{N} \rangle\right)=\langle f\Psi_{M}(m),\tilde{\mu}_{t} \rangle + \limsup_{N}\mathbb{E}\left(\langle f(1-\Psi_{M}(m)),\tilde{\mu}_{t}^{N} \rangle\right)
\end{equation}
and
\begin{equation}\liminf_{N}\mathbb{E}\left(\langle f,\tilde{\mu}_{t}^{N} \rangle\right)=
\langle f\Psi_{M}(m),\tilde{\mu}_{t} \rangle +
\liminf_{N}\mathbb{E}\left(\langle f(1-\Psi_{M}(m)),
\tilde{\mu}_{t}^{N} \rangle\right).
\end{equation}
Since $\mu_{t}^{N}=\tilde{\mu}_{t}^{N}\circ m^{-1}$ and $\varphi
\geq 1,$ we can write
\begin{align*}\left|\langle f(1-\Psi_{M}(m)),\tilde{\mu}_{t}^{N} \rangle\right|&=\left|\langle f(1-\Psi_{M}),\mu_{t}^{N} \rangle\right|\\
&\leq \|f\|_{\infty}\left|\langle \varphi(1-\Psi_{M}),\mu_{t}^{N}
\rangle\right|
\end{align*}
Hence, \begin{equation}\mathbb{E}\left(\left|\langle
f(1-\Psi_{M}(m)),\tilde{\mu}_{t}^{N} \rangle\right|\right) \leq
\|f\|_{\infty}\mathbb{E}\left(\left|\langle
\varphi(1-\Psi_{M}),\mu_{t}^{N} \rangle\right|\right).
\end{equation}
Now, $\left|\langle \varphi(1-\Psi_{M}),\mu_{t}^{N} \rangle\right|$
is bounded for all $N$. Indeed, since $\varphi$ is sublinear, for
all $t,$ we have
\[\langle \varphi,\mu_{t}^{N}\rangle \leq \langle \varphi,\mu_{0}^{N}\rangle.\]
Hence,
\[\langle \varphi(1-\Psi_{M}),\mu_{t}^{N}\rangle \leq \langle \varphi,\mu_{0}^{N}\rangle.\]
Now fix $\epsilon>0.$ Since $\langle \varphi,\mu_{0}^{N}\rangle \to
\langle \varphi,\mu_{0}\rangle$ as $N\to\infty$, we can find
$N_{0}>0$ such that for all $N\geq N_{0},$
\[\langle \varphi,\mu_{0}^{N}\rangle \leq
\langle \varphi,\mu_{0}\rangle+ \epsilon.\] Let $M_{0}=\max_{N\leq
N_{0}}\langle \varphi,\mu_{0}^{N}\rangle$ and set
$C_{0}=\max(M_{0},\langle \varphi,\mu_{0}\rangle+ \epsilon).$ For
all $N$ we have,
\[\langle \varphi,\mu_{t}^{N}\rangle \leq C_{0}\] and a fortiori,
\[\langle \varphi(1-\Psi_{M}),\mu_{t}^{N}\rangle \leq C_{0}\] so it
is bounded for all $N$. Moreover, it converges to $\left|\langle
\varphi(1-\Psi_{M}),\mu_{t} \rangle\right|$ in probability as $N \to
\infty.$ Thus, by the Bounded convergence Theorem,
\begin{align*}
\mathbb{E}\left(\left|\langle \varphi(1-\Psi_{M}),\mu_{t}^{N} \rangle\right|\right)&\to \mathbb{E}\left(\left|\langle \varphi(1-\Psi_{M}),\mu_{t}\rangle\right|\right)\\
&\qquad=\left|\langle \varphi(1-\Psi_{M}),\mu_{t}\rangle\right|
\end{align*}
as $N \to \infty.$ Taking the limsup over $N$ in the
expression (20) we obtain
\[\limsup_{N} \mathbb{E}\left(\left|\langle f(1-\Psi_{M}(m)),\tilde{\mu}_{t}^{N} \rangle\right|\right)\leq
\|f\|_{\infty}\left|\langle
\varphi(1-\Psi_{M}),\mu_{t}\rangle\right|.\]
As $M \to \infty$, $\varphi(1-\Psi_{M}) \searrow 0$ and is positive.
Hence, by the Monotone Convergence Theorem,
\[\left|\langle \varphi(1-\Psi_{M}),\mu_{t} \rangle\right| \to 0 \mbox{ as } M\to \infty.\]
So \[\lim_{M}\limsup_{N} \mathbb{E}\left(\left|\langle
f(1-\Psi_{M}(m)),\tilde{\mu}_{t}^{N}
\rangle\right|\right)=\lim_{M}\liminf_{N}
\mathbb{E}\left(\left|\langle f(1-\Psi_{M}(m)),\tilde{\mu}_{t}^{N}
\rangle\right|\right)=0\] Moreover,
\[\left|\langle f(1-\Psi_{M}(m)),\tilde{\mu}_{t} \rangle\right|\leq \|f\|_{\infty}\left|\langle \varphi(1-\Psi_{M}),\mu_{t} \rangle\right|\] that is going to 0 as $M \to \infty.$
Thus, \[\lim_{M \to \infty}\langle
f\Psi_{M}(m),\tilde{\mu}_{t}\rangle=\langle
f,\tilde{\mu}_{t}\rangle.\] Hence, taking the limit as $M \to
\infty$ in the relations (18) and (19) we obtain
\[\limsup_{N}\mathbb{E}\left(\langle f,\tilde{\mu}_{t}^{N} \rangle\right)=\liminf_{N}\mathbb{E}\left(\langle f,\tilde{\mu}_{t}^{N} \rangle\right)=\langle f,\tilde{\mu}_{t}\rangle.\]
Thus,
\[\lim_{N}\mathbb{E}\left(\langle f,\tilde{\mu}_{t}^{N} \rangle\right)=\langle f,\tilde{\mu}_{t}\rangle\] as required.
\begin{flushright}
$\square$
\end{flushright}
{\it Proof of Lemma 3.2 : }Set $f^{M}=f\Psi_{M}(m).$ We can write,
\begin{displaymath}
\mathbb{E}\left(\langle f^{M},\tilde{\mu}_{t}^{N} \rangle
\right)=\frac{2^{-q(\tau)}}{N} \sum_{(i_{1},\ldots, i_{n})\in[N]_{[n]}^{\star}} \mathbb{E}\bigg(f^{M}
\left(\xi_{t}^{I(i_{1},\ldots,i_{n})}\right)1_{\{S_{I(i_{1},\ldots,i_{n})}\leq
t < T_{I(i_{1},\ldots,i_{n})}\}}\bigg)
\end{displaymath}
where $I:[N]_{[n]}^{\star}\to \mathbb{T}_{n}^{\star}[N]$ is the map defined in subsection 2.1. For $(i_{1},\ldots, i_{n}) \in[N]_{[n]}^{\star},$ let $f_{I(i_{1},\ldots, i_{n})}^{M}$ be the corresponding element for $f^{M}$ in
$C_{b}(A_{I(i_{1},\ldots, i_{n})}(0,t))$ as defined in subsection 2.2. Then,
\begin{align*}
&\mathbb{E}\left(\langle f^{M},\tilde{\mu}_{t}^{N} \rangle
\right)\\
&=\frac{1}{N} \sum_{(i_{1},\ldots, i_{n}) \in[N]_{[n]}^{\star}} \mathbb{E}\bigg(f_{I(i_{1},\ldots, i_{n})}^{M}
\left(\xi_{t}^{I(i_{1},\ldots,i_{n})}\right)1_{\{S_{I(i_{1},\ldots,i_{n})}\leq
t < T_{I(i_{1},\ldots,i_{n})}\}}\bigg)\\
&=\frac{1}{N^{n}} \sum_{(i_{1},\ldots, i_{n}) \in[N]_{[n]}^{\star}} \mathbb{E}\Bigg[\int_{A_{I(i_{1},\ldots,i_{n})}^{y_{(i_{1},\ldots,i_{n})}}(0,t)}f_{I(i_{1},\ldots, i_{n})}^{M}\left(\xi\right)
K_{\xi}\\
&\exp\left(-\int_{\Delta(\xi)}\int_{E}K(y_{r},y^{\prime})\mu_{\Pi_{\xi}(r)}^{N,[N]\backslash\{i_{1},\ldots, i_{n}\}}(dy^{\prime})dr
\right)
\exp\left(-\int_{0}^{t}\frac{K_{s}(\xi)}{N}ds
\right)\nu_{I(i_{1},\ldots,i_{n})}^{y_{(i_{1},\ldots,i_{n})}}(d\xi)\Bigg]
\end{align*}
using the expression of the conditional expectation obtained at the end of the subsection 2.4.2.
\paragraph{Step 1:} We are going to give an alternative expression of $\mathbb{E}\left(\langle f^{M},\tilde{\mu}_{t}^{N} \rangle \right)$ using a new measure. Define on $(0,\infty)^{n}$ the measure,
\begin{displaymath}
\mu_{0}^{n,N}=\frac{1}{N^{n}}\sum_{(i_{1},\ldots,i_{n}) \in [N]_{[n]}^{\star}}\delta_{(y_{i_{1}},\ldots,y_{i_{n}})}.
\end{displaymath}
For $(i_{1},\ldots, i_{n}) \in[N]_{[n]}^{\star}$ define
\begin{align*}
&G(y_{i_{1}},\ldots,y_{i_{n}})\\
&=\mathbb{E}\Bigg[\int_{A_{I(i_{1},\ldots, i_{n})}^{y_{(i_{1},\ldots,i_{n})}}(0,t)}f_{I(i_{1},\ldots, i_{n})}^{M}\left(\xi\right)
K_{\xi}\exp\left(-\int_{0}^{t}\frac{K_{s}(\xi)}{N}ds
\right)\\
&\exp\left(-\int_{\Delta(\xi)}\int_{E}K(y_{r},y^{\prime})\mu_{\Pi_{\xi}(r)}^{N,[N]\backslash\{i_{1},\ldots, i_{n}\}}(dy^{\prime})dr
\right)
\nu_{I(i_{1},\ldots, i_{n})}^{y_{(i_{1},\ldots,i_{n})}}(d\xi)\Bigg].
\end{align*}
In this expression, the first three terms in the integral (inside the expectation) only depends on $(y_{i_{1}},\ldots,y_{i_{n}})$, whereas the last term depend on $(y_{1},\ldots,y_{N})$. Nevertheless we only need to know which particles are present in the first three terms to work out which particles are going to be in the last term. Moreover, if there exist $i_{1},i_{1}^{\prime}\in [N]$ with $i_{1}\neq i_{1}^{\prime}$ but $y_{i_{1}}=y_{i_{1}^{\prime}}$ then we have
\[G(y_{i_{1}^{\prime}},y_{i_{2}}\ldots,y_{i_{n}})=G(y_{i_{1}},y_{i_{2}}\ldots,y_{i_{n}}).\] Thus we can write,
\begin{align*}
\mathbb{E}\left(\langle f^{M},\tilde{\mu}_{t}^{N} \rangle
\right)&=\frac{1}{N^{n}} \sum_{(i_{1},\ldots, i_{n}) \in[N]_{[n]}^{\star}}G(y_{i_{1}},\ldots,y_{i_{n}})\\
&=\int_{(0,\infty)^{n}}G(x_{1},\ldots,x_{n})\mu_{0}^{n,N}(dx_{1}\ldots dx_{n}).
\end{align*}
\paragraph{Step 2:}
We are going to show that we can write
\begin{displaymath}
\mathbb{E}\left(\langle f^{M},\tilde{\mu}_{t}^{N}\rangle
\right)
=\int_{(0,\infty)^{n}}G(x_{1},\ldots,x_{n})
\mu_{0}^{N}(dx_{1})\ldots \mu_{0}^{N}(dx_{n})
+\circ(\frac{1}{N}).
\end{displaymath}
Let us consider
\begin{displaymath}
\left(\mu_{0}^{N}\right)^{\otimes
n}=\left(\frac{1}{N}\sum_{i=1}^{N}\delta_{y_{i}}\right)^{n}
=\mu_{0}^{n,N}+C_{0}^{n,N}
\end{displaymath}
where $C_{0}^{n,N}$ is a finite sum of terms of the form
\[\mu_{0}^{N,n_{1},\ldots,n_{k}}=\frac{1}{N^{n}}\sum_{\substack{i_{1},\ldots,i_{k}=1 \\ $
\scriptsize{with} $ i_{1},\ldots,i_{k} \\ $ \scriptsize{distinct}
$}}^{N} \delta_{(\underbrace{y_{i_{1}},
\ldots,y_{i_{1}}}_{n_{1}~times},\underbrace{y_{i_{2}},\ldots,y_{i_{2}}}_{n_{2}~times},\ldots,
\underbrace{y_{i_{k}},\ldots,y_{i_{k}}}_{n_{k}~times})}\] for $k>1$
with $n_{1}+n_{2}+\ldots+n_{k}=n$ modulo some permutations of the masses. Take
$k<n-1$ and consider,
\begin{align*}
a_{N}^{k}=\frac{1}{N^{n}}\sum_{\substack{i_{1},\ldots,i_{k}=1 \\ $
\scriptsize{with} $ i_{1},\ldots,i_{k} \\ $ \scriptsize{distinct}
$}}^{N}&\mathbb{E}\bigg(f^{M}
\left(\xi_{t}^{I(i_{1},
\ldots,i_{1},\ldots,i_{k},\ldots,i_{k})} \right)\\
&\qquad\quad 1_{\{S_{I(i_{1},
\ldots,i_{1},\ldots,i_{k},\ldots,i_{k})}\leq t <
T_{I(i_{1}, \ldots,i_{1},\ldots,
i_{k},\ldots,i_{k})}\}}\bigg)N^{n-1}.
\end{align*}
where $I(i_{1}, \ldots,i_{1},\ldots,
i_{k},\ldots,i_{k})$ is obtained from $I(i_{1},\ldots,i_{n})$ with $(i_{1},\ldots,i_{n})\in [N]_{[n]}^{\star}$ by substituting in the tree $i(i_{1},\ldots,i_{n})$, the particles $\{i_{p}:1\leq p \leq n_{1}\}$ by the particle $i_{1}$, the particles $\{i_{p}:n_{1}+1\leq p \leq n_{1}+n_{2}\}$ by the particle $i_{2}$\ldots, and finally the particles $\{i_{p}:n-n_{k}+1\leq p \leq n\}$ by the particle $i_{k}$. It is a tree where particles are not all distinct. We are going to prove that this quantity is $\circ(\frac{1}{N}).$
By a similar argument to before we can write
\[a_{N}^{k}=\frac{1}{N^{n}}\sum_{\substack{i_{1},\ldots,i_{k}=1 \\ $
\scriptsize{with} $ i_{1},\ldots,i_{k} \\ $ \scriptsize{distinct}
$}}^{N}g(y_{i_{1}},\ldots y_{i_{k}})\]
where $g:(0,\infty)^{k}\to\mathbb{R}$ is the expectation above.

Now define
\[\left(\mu_{0}^{N}\right)^{\otimes
k}=\left(\frac{1}{N}\sum_{i=1}^{N}\delta_{y_{i}}\right)^{k}\] and on
$(0,\infty)^{k}$ the measure,
\begin{displaymath}
\mu_{0}^{k,N}=\frac{1}{N^{k}}\sum_{\substack{i_{1},\ldots,i_{k}=1 \\
$ \scriptsize{with} $ i_{1},\ldots,i_{k} \\ $ \scriptsize{distinct}
$}}^{N}\delta_{(y_{i_{1}},\ldots,y_{i_{k}})}.
\end{displaymath}
Then,
\[\mu_{0}^{k,N}\leq \left(\mu_{0}^{N}\right)^{\otimes
k}.\] We can identify $\mu_{0}^{N,n_{1},\ldots,n_{k}}$ and
$\mu_{0}^{k,N}$. Hence,
\begin{align*}
a_{N}^{k}&=\frac{1}{N^{n-k}}\int_{(0,\infty)^{k}}g(x_{1},\ldots
x_{k})\mu_{0}^{k,N}(dx_{1}\ldots dx_{k})\\
&\leq \frac{1}{N^{n-k}}\int_{(0,\infty)^{k}}g(x_{1},\ldots
x_{k})\mu_{0}^{N}(dx_{1})\ldots \mu_{0}^{N}(dx_{k})
\end{align*}
For $i\in\mathbb{T}_{n}^{\star}[n]$, let $A_{i,M}(0,t)=\{\xi \in A_{i}(0,t):
m(\xi)\leq M+1\}$. On $A_{i,M}(0,t)$, the map $A(0,t) \to
(0,\infty): \xi \to K_{\xi}$ is bounded by
$C=\left(\sup_{y,y^{\prime} \in
[0,M]^{2}}K(y,y^{\prime})\right)^{n-1}$ ( which is attained as $K$
is bounded on $[0,M]^{2}$ compact). Outside
$A_{i,M}(0,t)$, $\Psi_{M}(m)$ is the zero-function and so
is $f_{M}.$
Now, for all $(i_{1},\ldots,i_{n})\in [N]^{n}$ non necessarily distinct,
\[
g(y_{i_{1}},\ldots,y_{i_{k}})
=\mathbb{E}\left(P_{t}^{N}f_{I(i_{1},\ldots,i_{n})}^{M}\left(\tilde{y}_{(i_{1},\ldots,i_{n})}\right)\right).
\]
Now,
\begin{align*}
&\left|P_{t}^{N}f_{I(i_{1},\ldots,i_{n})}^{M}\left(\tilde{y}_{(i_{1},\ldots,i_{n})}\right)\right|
\\
&=\Big|\int_{A_{I(i_{1},\ldots,i_{n})}^{y_{(i_{1},\ldots,i_{n})}}(0,t)}f_{I(i_{1},\ldots,i_{n})}^{M}\left(\xi\right)
K_{\xi}\exp\left(-\int_{0}^{t}\frac{K_{s}(\xi)}{N}ds \right)\\
&\quad\exp\left(-\int_{\Delta(\xi)}\int_{E}K(y_{r},y^{\prime})\mu_{\Pi_{\xi}(r)}^{N,J}(dy^{\prime})dr
\right)
\nu_{I(i_{1},\ldots,i_{n})}^{y_{(i_{1},\ldots,i_{n})}}(d\xi)\Big|\\
&\leq C
\|f\|_{\infty}\int_{A_{I(i_{1},\ldots,i_{n})}^{y_{(i_{1},\ldots,i_{n})}}(0,t)}
\nu_{I(i_{1},\ldots,i_{n})}^{y_{(i_{1},\ldots,i_{n})}}(d\xi)\\
&=C \|f\|_{\infty}\int_{(0,\infty)^{n}}
\sum_{\substack{$\scriptsize{ Possible permutations }$\\
$\scriptsize{ for }$ s_{1},\ldots,s_{n-1}}}1_{\{s_{1}\leq \cdots
\leq s_{n-1}\}}ds_{1}\ldots ds_{n-1}\\
&\leq C \|f\|_{\infty}t^{n-1}.
\end{align*}
Hence,
\begin{align*}
a_{N}^{k}&\leq\frac{{C\|f
\|}_{\infty}t^{n-1}}{N^{n-k}} \int_{(0,\infty)^k}\mu_{0}^{N}(dx_{1})\ldots \mu_{0}^{N}(dx_{k})\\
&\leq\frac{{C\|f
\|}_{\infty}t^{n-1}}{N^{n-k}}
\int_{(0,\infty)^k}\varphi(x_{1})\mu_{0}^{N}(dx_{1})\ldots \varphi(x_{n})\mu_{0}^{N}(dx_{k})\\
&=\frac{{C\|f\|}_{\infty}t^{n-1}}{N^{n-k}} \langle
\varphi,\mu_{0}^{N}\rangle ^{k}.
\end{align*}
as $\varphi \geq 1$. It is clear that this quantity is going to $0$
as $N\to\infty$ since $k<n-1$ and $\langle
\varphi,\mu_{0}^{N}\rangle \to \langle \varphi,\mu_{0}\rangle <
\infty$ as $N \to \infty$. Hence,
\begin{displaymath}
\mathbb{E}\left(\langle f_{M},\tilde{\mu}_{t}^{N}\rangle
\right)
=\int_{(0,\infty)^{n}}G(x_{1},\ldots,x_{n})
\mu_{0}^{N}(dx_{1})\ldots \mu_{0}^{N}(dx_{n})
+\circ(\frac{1}{N}).
\end{displaymath}
as required.
\paragraph{Step 3:}
Let $i_{0} \in \mathbb{T}_{n}^{\star}[n]$ with type $\tau$ and without associated masses. For each $(x_{1},\ldots,x_{n})\in (0,\infty)^{n}$, we can associate to $i_{0}$ the vector of masses $(x_{1},\ldots,x_{n})$. Hence we can write
\[G(x_{1},\ldots,x_{n})=\mathbb{E}\left(P_{t}^{N}f_{i_{0}}^{M}\left(x_{1},\ldots,x_{n}\right)\right)\]
By Theorem 2.4, for all $x_{1},\ldots,x_{n}\in (0,\infty)$,
\[P_{t}^{N}f_{i_{0}}^{M}\left(x_{1},\ldots,x_{n}\right)\to
P_{t}f_{i_{0}}^{M}\left(x_{1},\ldots,x_{n}\right)\] in
probability as $N\to\infty$. Also,
$P_{t}^{N}f_{i_{0}}^{M}(x_{1},\ldots,x_{n})$ is
bounded by $C\|f\|_{\infty}t^{n-1}$. Hence, by the Bounded
Convergence Theorem,
\[\mathbb{E}\left(\langle f_{M},\tilde{\mu}_{t}^{N}\rangle\right)\to \int_{(0,\infty)^{n}}
P_{t}f_{i_{0}}^{M}\left((x_{1},\ldots,x_{n})\right)
\mu_{0}(dx_{1})\ldots \mu_{0}(dx_{n})\] as $N\to\infty$.

Now, we need to prove that
\[\int_{(0,\infty)^{n}}P_{t}f_{i_{0}}^{M}\left((x_{1},\ldots,x_{n})\right)
\mu_{0}(dx_{1})\ldots\mu_{0}(dx_{n})=\langle
f^{M},\tilde{\mu}_{t}\rangle.\] But,
\begin{align*}
&\int_{(0,\infty)^{n}}P_{t}f_{i_{0}}^{M}\left(x_{1},\ldots,x_{n}\right)
\mu_{0}(dx_{1})\ldots \mu_{0}(dx_{n})\\
&\qquad=\int_{(0,\infty)^{n}}\int_{A_{i_{0}}^{(x_{1},\ldots,x_{n})}(0,t)}K_{\xi}f_{i_{0}}^{M}\left(\xi\right)
\nu_{i_{0}}^{(x_{1},\ldots,x_{n})}(d\xi)\mu_{0}(dx_{1})\ldots \mu_{0}(dx_{n})\\
&\qquad\qquad\exp\left(-\int_{\Delta\left(\xi\right)}\int_{E}K(y_{r},y^{\prime})\mu_{r}(dy^{\prime})dr\right)\\
&\qquad=\int_{A_{i_{0}}(0,t)}f_{i_{0}}^{M}\left(\xi\right)\tilde{\mu}_{t}^{\prime}(d\xi)
=\int_{A_{\tau}(0,t)}f^{M}\left(\xi\right)\tilde{\mu}_{t}(d\xi)
=\langle f^{M}, \tilde{\mu}_{t} \rangle
\end{align*}
using Appendix 5.3. Hence,
\[\mathbb{E}(\langle f^{M},\tilde{\mu}_{t}^{N}\rangle)\to \langle f^{M},\tilde{\mu}_{t}\rangle \] as $N\to\infty$ as required.
\begin{flushright}
$\square$
\end{flushright}

\subsection{Convergence of the expectation of the square} We are still
working with $[N]=\{1,\ldots,N\}$ with associated masses
$y_{1},\ldots,y_{N}>0$. Take $\tau \in \mathbb{T}$ with $n$ leaves
that is $n(\tau)=n$ for some $n>0$. Fix $f\in C_{b}(A_{\tau}(0,t))$.
The aim in this section is to prove the following result.
\begin{proposition}
\[\mathbb{E}\left(\langle f,\tilde{\mu}_{t}^{N} \rangle ^{2}\right)\to \langle f,\tilde{\mu}_{t}\rangle^{2}\] as $N \to \infty$.
\end{proposition}

In this aim, we are going to proceed similarly to sections 2 and 3.
First we are going to express $\mathbb{E}\left(\langle
f,\tilde{\mu}_{t}^{N} \rangle ^{2}\right)$ as the expected value of
a finite sum of conditional expectations, then we will compute these
conditional expectations and give some properties about their
convergence. Finally we will prove Proposition 3.3.

\subsubsection{An expression for the expectation}

 Take $i \in
\mathbb{T}_{n}^{\star}[n]$ with type $\tau$ and associate the
function $I:[N]_{[n]}^{\star}\to \mathbb{T}_{n}^{\star}[N]$ defined
in 2.1. The aim is to compute 
\[
\mathbb{E}\left(\langle f,\tilde{\mu}_{t}^{N} \rangle ^{2}\right).
\]
We can write,
\begin{align*}
&\mathbb{E}\left(\langle
f,\tilde{\mu}_{t}^{N}\rangle^{2}\right)
=\mathbb{E}\left(\left(\frac{2^{-q(\tau)}}{N}\sum_{\sigma\in S_{n,N}}
f\left(\xi_{t}^{\sigma(i)}\right)1_{\{S_{\sigma(i)}
\leq t <
T_{\sigma(i)}\}}\right)^{2}\right)\\
&=\frac{2^{-2q(\tau)}}{N^{2}}\sum_{\sigma\in S_{n,N}}
\mathbb{E}\bigg(f^{2}\left(\xi_{t}^{\sigma(i)}\right)1_{\{S_{\sigma(i)}\leq
t < T_{\sigma(i)}\}}\bigg)\\
&+\frac{2^{-2q(\tau)}}{N^{2}}\sum_{\substack{\sigma,\sigma^{\prime}\in S_{n,N}\\ $\scriptsize{with }$ Im(\sigma)\cap Im(\sigma^{\prime})=\emptyset}}
\mathbb{E}\bigg(f\left(\xi_{t}^{\sigma(i)}\right)f\left(\xi_{t}^{\sigma^{\prime}(i)}\right)
1_{\{S_{\sigma(i)}\leq
t <
T_{\sigma(i)}\}}1_{\{S_{\sigma^{\prime}(i)}\leq
t < T_{\sigma^{\prime}(i)}\}}\bigg).
\end{align*}
where $Im(\sigma)=\{\sigma(j): j \in [N]\}.$
Observe that in the formula above we omitted terms of the form
\[
\mathbb{E}\bigg(f\left(\xi_{t}^{\sigma(i)}\right)f\left(\xi_{t}^{\sigma^{\prime}(i)}\right)
1_{\{S_{\sigma(i)}\leq t <
T_{\sigma(i)}\}}1_{\{S_{\sigma^{\prime}(i)}\leq
t < T_{\sigma^{\prime}(i)}\}}\bigg)
\]
for $\sigma, \sigma^{\prime} \in S_{n,N}$ with $Im(\sigma)\cap Im(\sigma^{\prime}) \neq
\emptyset$. Indeed, in this case, by the way we constructed the
$\{S_{j}: j \in \mathbb{T}_{+}^{\star}[N]\}$, the events
$\{S_{\sigma(i)}\leq t <
T_{\sigma(i)}\}$ and
$\{S_{\sigma^{\prime}(i)}\leq t <
T_{\sigma^{\prime}(i)}\}$ are disjoint (because by
construction the trees have all distinct leaves) and so the quantity
above is equal to $0$. Let us look at the first term in the
expression of $\mathbb{E}\left(\langle f,\tilde{\mu}_{t}^{N} \rangle
^{2}\right)$. We have
\[\frac{2^{-2q(\tau)}}{N^{2}}\sum_{\sigma\in S_{n,N}}
\mathbb{E}\bigg(f^{2}\left(\xi_{t}^{\sigma(i)}\right)1_{\{S_{\sigma(i)}\leq
t < T_{\sigma(i)}\}}\bigg)=
\frac{2^{-q(\tau)}}{N}\mathbb{E}\left(\langle
f^{2},\tilde{\mu}_{t}^{N} \rangle \right).\] From subsection $3.1$,
as $N \to \infty,$
\[\mathbb{E}\left(\langle f^{2},\tilde{\mu}_{t}^{N} \rangle \right) \to \langle f^{2},\tilde{\mu}_{t}\rangle <\infty.\]

Hence as $N \to \infty$,
\[\frac{2^{-2q(\tau)}}{N^{2}}\sum_{\sigma\in S_{n,N}}
\mathbb{E}\bigg(f^{2}\left(\xi_{t}^{\sigma(i)}\right)1_{\{S_{\sigma(i)}\leq
t < T_{\sigma(i)}\}}\bigg) \to 0.\] For
$\sigma,\sigma^{\prime}\in S_{n,N}$ with $Im(\sigma)\cap Im(\sigma^{\prime})=
\emptyset$, set
\[
J(\sigma(i),\sigma^{\prime}(i))=[N]\setminus
\{\lambda\left(\sigma(i)\right),\lambda\left(\sigma^{\prime}(i)\right)\}.
\]
Hence,
\begin{align*}
\mathbb{E}\left(\langle f,\tilde{\mu}_{t}^{N}\rangle^{2}\right)
&=\mathbb{E}\Bigg(\frac{2^{-2q(\tau)}}{N^{2}}\sum_{\substack{\sigma,\sigma^{\prime}\in S_{n,N}\\ $\scriptsize{with }$ Im(\sigma)\cap Im(\sigma^{\prime})=\emptyset}}
\mathbb{E}\bigg(f\left(\xi_{t}^{\sigma(i)}\right)f\left(\xi_{t}^{\sigma^{\prime}(i)}\right)\\
&\qquad 1_{\{S_{\sigma(i)}\leq t<T_{\sigma(i)}\}}
 1_{\{S_{\sigma^{\prime}(i)}\leq t < T_{\sigma^{\prime}(i)}\}}\vert
\mathcal{F}^{J(\sigma(i),\sigma^{\prime}(i))}
\bigg)\Bigg)+ \circ(\frac{1}{N})
\end{align*}
Thus, we need to compute the following conditional expectation
\begin{displaymath}
\mathbb{E}\bigg(f\left(\xi_{t}^{\sigma(i)}\right)f\left(\xi_{t}^{\sigma^{\prime}(i)}\right)
1_{\{S_{\sigma(i)}\leq t <
T_{\sigma(i)}\}}
1_{\{S_{\sigma^{\prime}(i)}\leq t < T_{\sigma^{\prime}(i)}\}}\vert
\mathcal{F}^{J(\sigma(i),\sigma^{\prime}(i))}
\bigg)
\end{displaymath}
for $\sigma,\sigma^{\prime} \in S_{n,N}$ with $Im(\sigma)\cap Im(\sigma^{\prime})=\emptyset$.
\subsubsection{An expression for the conditional expectation}
In order to compute $\mathbb{E}\left(\langle
f,\tilde{\mu}_{t}^{N}\rangle^{2}\right),$ without loss of
generality, by subsection 3.2.1, we need to compute
\[\mathbb{E}_{y^{1},y^{2},\mu_{0}^{N,J}}\left(f\left(\xi_{t}^{i_{1}}\right)f\left(\xi_{t}^{i_{2}}\right)1_{\{S_{i_{1}}\leq t < T_{i_{1}}\}}1_{\{S_{i_{2}}\leq t < T_{i_{2}}\}}|\mathcal{F}_{J} \right)\]
for $i_{1} \in \mathbb{T}_{n}^{\star}[n],
i_{2} \in \mathbb{T}_{n}^{\star}[2n]\backslash[n]$ both
with type $\tau$. The quantities $y^{1}$ and $y^{2}$ represent the
respective vector of masses (as defined in 2.2) for $i_{1}$ and
$i_{2}$, $J=\{2n+1,\ldots,N\}$ and
$\mathcal{F}_{J}=\sigma\left(U_{j}: j\in
\mathbb{T}_{+}^{\star}J\right).$ Define $f_{i_{1}}=f\circ
g_{i_{1}}$ and $f_{i_{2}}=f\circ
g_{i_{2}}$ where $g_{i_{1}}$ and
$g_{i_{2}}$ are the map on forgetting labels defined in subsection
2.2. Then $f_{i_{1}} \in
C_{b}\left(A_{i_{1}}(0,t)\right)$ and  $f_{i_{2}} \in
C_{b}\left(A_{i_{2}}(0,t)\right)$.
\paragraph{Step 1 : A symmetry argument} ~~\\~~\\
By a similar symmetry argument that the one we used in subsection 2.2, we can
write
\begin{align*}
&\mathbb{E}_{y^{1},y^{2},\mu_{0}^{N,J}}\bigg(f\left(\xi_{t}^{i_{1}}\right)f\left(\xi_{t}^{i_{2}}\right)1_{\{S_{i_{1}}\leq t < T_{i_{1}}\}}1_{\{S_{i_{2}}\leq t < T_{i_{2}}\}}|\mathcal{F}_{J} \bigg)\\
&\qquad=\mathbb{E}_{y^{1},y^{2},\mu_{0}^{N,J}}\bigg(f_{i_{1}}\left(\xi_{t}^{i_{1}}\right)f_{i_{2}}\left(\xi_{t}^{i_{2}}\right)1_{\{S_{i_{1}}\leq
t < T_{i_{1}}\}}1_{\{S_{i_{2}}\leq t <
T_{i_{2}}\}}|\mathcal{F}_{J} \bigg).
\end{align*}
Thus it is enough to compute
\[\mathbb{E}_{y^{1},y^{2},\mu_{0}^{N,J}}\left(f_{i_{1}}\left(\xi_{t}^{i_{1}}\right)f_{i_{2}}\left(\xi_{t}^{i_{2}}\right)1_{\{S_{i_{1}}\leq t < T_{i_{1}}\}}1_{\{S_{i_{2}}\leq t < T_{i_{2}}\}}|\mathcal{F}_{J} \right).\]
\paragraph{Step 2 : Simplifying the computation}~~\\~~\\
Set $\mathcal{F}_{[n]}=\sigma\left(U_{j}:j \in
\mathbb{T}_{+}^{\star}[n]\right)$ and $\mathcal{F}_{[2n]\backslash
[n]}=\sigma\left(U_{j}:j \in \mathbb{T}_{+}^{\star}[2n]\backslash
[n]\right)$ and look at
\begin{equation}\mathbb{E}_{y^{1},y^{2},\mu_{0}^{N,J}}\left(f_{i_{1}}\left(\xi_{t}^{i_{1}}\right)f_{i_{2}}\left(\xi_{t}^{i_{2}}\right)1_{\{S_{i_{1}}\leq t < T_{i_{1}}\}}1_{\{S_{i_{2}}\leq t < T_{i_{2}}\}}|\mathcal{F}_{J}\vee \mathcal{F}_{[n]} \vee \mathcal{F}_{[2n]\backslash [n]} \right).
\end{equation}
Observe that
$f_{i_{1}}\left(\xi_{t}^{i_{1}}\right)1_{\{S_{i_{1}}\leq
t < T_{i_{1}}^{[n]}\}}$ is measurable with respect to
$\mathcal{F}_{[n]}$ and that
\\$f_{i_{2}}\left(\xi_{t}^{i_{2}}\right)1_{\{S_{i_{2}}\leq
t < T_{i_{2}}^{[2n]\backslash [n]}\}}$ is measurable with
respect to $\mathcal{F}_{[2n]\backslash [n]}.$ Moreover, $\{t <
T_{i_{1}}\} \subset \{t < T_{i_{1}}^{[n]}\}$
and $\{t < T_{i_{2}}\} \subset \{t <
T_{i_{2}}^{[2n]\backslash [n]}\}$. Thus, the quantity
(21) is equal to
\begin{align*}
&f_{i_{1}}\left(\xi_{t}^{i_{1}}\right)f_{i_{2}}\left(\xi_{t}^{i_{2}}\right)1_{\{S_{i_{1}}\leq t < T_{i_{1}}^{[n]}\}}1_{\{S_{i_{2}}\leq t < T_{i_{2}}^{[2n]\backslash [n]}\}}\\
&\qquad\mathbb{E}_{y^{1},y^{2},\mu_{0}^{N,J}}\left(1_{\{t <
T_{i_{1}}\}}1_{\{t <
T_{i_{2}}\}}|\mathcal{F}_{J}\vee \mathcal{F}_{[n]} \vee
\mathcal{F}_{[2n]\backslash [n]} \right).
\end{align*}
Hence we can write,
\begin{align*}
&\mathbb{E}_{y^{1},y^{2},\mu_{0}^{N,J}}\left(f_{i_{1}}\left(\xi_{t}^{i_{1}}\right)f_{i_{2}}\left(\xi_{t}^{i_{2}}\right)1_{\{S_{i_{1}}\leq t < T_{i_{1}}\}}1_{\{S_{i_{2}}\leq t < T_{i_{2}}\}}|\mathcal{F}_{J}\right)\\
&\qquad= \int_{A_{i_{1}}^{y^{1}}(0,t)\times A_{i_{1}}^{y^{2}}(0,t)}f_{i_{1}}\left(\xi_{1}\right)f_{i_{2}}\left(\xi_{2}\right)\\
&\qquad\qquad\mathbb{P}_{y^{1},y^{2},\mu_{0}^{N,J}}\Big(\{t <
T_{i_{1}}\}\cap \{t <
T_{i_{2}}\}|\mathcal{F}_{J}\vee \mathcal{F}_{[n]} \vee
\mathcal{F}_{[2n]\backslash [n]}\Big)\left(\xi_{1},\xi_{2}\right)
\\
&\qquad\qquad\mathbb{P}_{y^{1},y^{2}}\left(S_{i_{1}}\leq
t < T_{i_{1}}^{[n]},S_{i_{2}}\leq t <
T_{i_{2}}^{[2n]\backslash [n]},\xi_{t}^{i_{1}}
\in d\xi_{1},\xi_{t}^{i_{2}} \in d\xi_{2}\right).
\end{align*}

\begin{lemma}
\begin{enumerate}
\item Conditional on the event $\{S_{i_{1}}\leq t < T_{i_{1}}^{[n]}\}\cap \{S_{i_{2}}\leq t < T_{i_{2}}^{[2n]\backslash [n]}\}$,
\begin{align*}
&\mathbb{P}_{y^{1},y^{2},\mu_{0}^{N,J}}\Big(\{t < T_{i_{1}}\}\cap \{t < T_{i_{2}}\}|\mathcal{F}_{J}\vee \mathcal{F}_{[n]} \vee \mathcal{F}_{[2n]\backslash [n]}\Big)\\
&\qquad= \exp\left(-\frac{1}{2}\int_{0}^{t}\frac{K_{s}(\xi_{t}^{i_{1}},\xi_{t}^{i_{2}})}{N} ds\right)\\
&\qquad\exp\left(-\int_{\Delta(\xi_{t}^{i_{1}})}\int_{E}K(y_{r},y^{\prime})\mu_{\Pi_{\xi_{t}^{i_{1}}}(r)}^{N,J}
(dy^{\prime})dr\right)\\
&\qquad
\exp\left(-\int_{\Delta(\xi_{t}^{i_{2}})}\int_{E}K(y_{r},y^{\prime})\mu_{\Pi_{\xi_{t}^{i_{2}}}(r)}^{N,J}(dy^{\prime})dr\right)
\end{align*}
where
\[K_{s}(\xi_{t}^{i_{1}},\xi_{t}^{i_{2}})=-\frac{1}{2}\sum_{\substack{r
\in \Pi_{\xi_{1}}^{-1}(s) \\ r^{\prime} \in \Pi_{\xi_{2}}^{-1}(s)
}}K(y_{r},y_{r^{\prime}}).\]
\item \begin{align*}
&\mathbb{P}_{y^{1},y^{2}}\left(S_{i_{1}}\leq t < T_{i_{1}}^{[n]},S_{i_{2}}\leq t < T_{i_{2}}^{[2n]\backslash [n]},\xi_{t}^{i_{1}} \in d\xi_{1},\xi_{t}^{i_{2}} \in d\xi_{2}\right)\\
&\qquad=\frac{K_{\xi_{1}}K_{\xi_{2}}}{N^{2n-2}}\exp\left(-\int_{0}^{t}\frac{K_{s}(\xi_{1})}{N}ds\right)
\exp\left(-\int_{0}^{t}\frac{K_{s}(\xi_{2})}{N}ds\right)\nu_{i_{1}}^{y^{1}}(d\xi_{1})\nu_{i_{2}}^{y^{2}}(d\xi_{2}).
\end{align*}
\end{enumerate}
\end{lemma}
In the lemma above, the quantity (1) is similar to the one in Lemma 2.2 except that
here, it means that for $\xi_{1}\in
A_{i_{1}}^{y^{1}}(0,t),\xi_{2}\in
A_{i_{2}}^{y^{2}}(0,t)$ and
$\omega^{J}=(\omega_{i}^{J}:i\in \mathbb{T}_{+}^{\star}J)$, the
configurations $\xi_{1}$ and $\xi_{2}$ do not interfere with
each other and are not killed by $\left(\mu_{r}^{N,J}\right)_{r \geq
0}(\omega^{J})$. Hence, mimicking the proof of Lemma 2.2, we obtain the
formula of the lemma. The second part of the proof comes from applying the Lemma 2.3 along with an argument of independence.

Hence,
\begin{align*}
&\mathbb{E}_{y^{1},y^{2},\mu_{0}^{N,J}}\left(f_{i_{1}}\left(\xi_{t}^{i_{1}}\right)f_{i_{2}}
\left(\xi_{t}^{i_{2}}\right)1_{\{S_{i_{1}}\leq t < T_{i_{1}}\}}1_{\{S_{i_{2}}\leq t < T_{i_{2}}\}}|\mathcal{F}_{J}\right)\\
&\qquad=\int_{A_{i_{1}}^{y^{1}}(0,t)\times A_{i_{2}}^{y^{2}}(0,t)}f_{i_{1}}\left(\xi_{1}\right)f_{i_{2}}\left(\xi_{2}\right)
\exp\left(-\frac{1}{2}\int_{0}^{t}\frac{K_{s}(\xi_{1},\xi_{2})}{N} ds\right)\\
&\qquad\qquad\frac{K_{\xi_{1}}K_{\xi_{2}}}{N^{2n-2}}\exp\left(-\int_{0}^{t}\frac{K_{s}(\xi_{1})}{N}ds\right)
\exp\left(-\int_{0}^{t}\frac{K_{s}(\xi_{2})}{N}ds\right)\\
&\qquad\qquad\exp\left(-\int_{\Delta(\xi_{1})}\int_{E}K(y_{r},y^{\prime})\mu_{\Pi_{\xi_{1}}(r)}^{N,J}(dy^{\prime})dr\right)\\
&\qquad\qquad
\exp\left(-\int_{\Delta(\xi_{2})}\int_{E}K(y_{r},y^{\prime})\mu_{\Pi_{\xi_{2}}(r)}^{N,J}(dy^{\prime})dr\right)\nu_{i_{1}}^{y^{1}}(d\xi_{1})
\nu_{i_{2}}^{y^{2}}(d\xi_{2}).
\end{align*}

\subsubsection{Convergence for the conditional expectation}
We keep the same notations that in the previous section and we define,
\[P_{t}^{N}f_{i_{1}}(y^{1})f_{i_{2}}(y^{2})=\mathbb{E}_{y^{1},y^{2},\mu_{0}^{N,J}}\left(f_{i_{1}}\left(\xi_{t}^{i_{1}}\right)f_{i_{2}}\left(\xi_{t}^{i_{2}}\right)1_{\{S_{i_{1}}\leq t < T_{i_{1}}\}}1_{\{S_{i_{2}}\leq t < T_{i_{2}}\}}|\mathcal{F}_{J} \right)N^{2n-2}.\]
and
\begin{align*}
&P_{t}f_{i_{1}}(y^{1})f_{i_{2}}(y^{2})\\
&=\int_{A_{i_{1}}^{y^{1}}(0,t)\times
A_{i_{2}}^{y^{2}}(0,t)}
K_{\xi_{1}}f_{i_{1}}(\xi_{1})\exp\left(-\int_{\Delta(\xi_{1})}\int_{E}K(y_{r},y^{\prime})\mu_{\Pi_{\xi_{1}}(r)}^{N,J}(dy^{\prime})dr\right)
\nu_{i_{1}}^{y^{1}}(d\xi_{1})\\
&\qquad
K_{\xi_{2}}f_{i_{2}}(\xi_{2})\exp\left(-\int_{\Delta(\xi_{2})}\int_{E}K(y_{r},y^{\prime})\mu_{\Pi_{\xi_{2}}(r)}^{N,J}(dy^{\prime})dr\right)
\nu_{i_{2}}^{y^{2}}(d\xi_{2}).
\end{align*}

\begin{proposition}
\begin{enumerate}
\item $P_{t}^{N}f_{i_{1}}(y^{1})f_{i_{2}}(y^{2})\to P_{t}f_{i_{1}}(y^{1})f_{i_{2}}(y^{2})$ as $N\to\infty$ in probability.
\item$P_{t}f_{i_{1}}(y^{1})f_{i_{2}}(y^{2})=P_{t}f_{i_{1}}(y^{1})P_{t}f_{i_{2}}(y^{2}).$
\end{enumerate}
\end{proposition}

We are not proving this proposition since the proof is similar to
the proof of Theorem 2.4.
\subsubsection{Proof of Proposition 3.3}
We are going to sketch the proof of this proposition as it is very
similar to the proof of Proposition 3.1.
\paragraph{Step 1:} For $M>0$, define $\Psi_{M}: (0,\infty) \to (0,\infty)$ to be
\begin{displaymath}
\Psi_{M}(x)=\left\{\begin{array}{ll} 1& \textrm{if $x \leq M$}\\
-x+(M+1)& \textrm{if $ M \leq x \leq M+1$}
\\ 0&
\textrm{if $x \geq M+1$}\end{array} \right.
\end{displaymath}
We can write,
\[\langle f,\tilde{\mu}_{t}^{N} \rangle=\langle f\Psi_{M}(m),\tilde{\mu}_{t}^{N} \rangle + \langle f(1-\Psi_{M}(m)),\tilde{\mu}_{t}^{N} \rangle.\]
Hence,
\[\langle f,\tilde{\mu}_{t}^{N} \rangle^{2}=\langle f\Psi_{M}(m),\tilde{\mu}_{t}^{N} \rangle^{2} + \langle f(1-\Psi_{M}(m)),\tilde{\mu}_{t}^{N} \rangle^{2}+2\langle f\Psi_{M}(m),\tilde{\mu}_{t}^{N} \rangle \langle f(1-\Psi_{M}(m)),\tilde{\mu}_{t}^{N} \rangle\] and so
\begin{align*}\label{key}
\mathbb{E}\left(\langle f,\tilde{\mu}_{t}^{N}
\rangle^{2}\right)&=\mathbb{E}\left(\langle
f\Psi_{M}(m),\tilde{\mu}_{t}^{N} \rangle^{2}\right) +
\mathbb{E}\left(\langle f(1-\Psi_{M}(m)),\tilde{\mu}_{t}^{N}
\rangle^{2}\right)\\
&\qquad+2\mathbb{E}\left(\langle f\Psi_{M}(m),\tilde{\mu}_{t}^{N}
\rangle \langle f(1-\Psi_{M}(m)),\tilde{\mu}_{t}^{N} \rangle\right).
\end{align*}
\paragraph{Step 2:} We prove that $\mathbb{E}\left(\langle f\Psi_{M}(m),\tilde{\mu}_{t}^{N} \rangle^{2}\right)\to\langle f\Psi_{M}(m),\tilde{\mu}_{t}\rangle^{2}$ as $N \to \infty.$
Let $i_{1}\in \mathbb{T}_{n}^{\star}[n]$ and
$i_{2}\in \mathbb{T}_{n}^{\star}[2n]\backslash[n]$ both
with type $\tau$, and associate the maps $I_{1} :
[N]_{[n]}^{\star}\to \mathbb{T}_{n}^{\star}[N]$ and $I_{2} :
[N]_{[2n]\backslash[n]}^{\star}\to \mathbb{T}_{n}^{\star}[N]$ defined in subsection 2.1. For $M>0$, let
$f^{M}=f\Psi_{M}(m)$ where $m$ is the mass function. Let
$f_{i_{1}}^{M}$ and $f_{i_{2}}^{M}$ be the correspondent
representant for $f_{M}$ respectively in
$C_{b}(A_{i_{1}}(0,t))$ and
$C_{b}(A_{i_{2}}(0,t))$. Define on $(0,\infty)^{2n}$, the
measure
\begin{displaymath}
\mu_{0}^{2n,N}=\frac{1}{N^{2n}}\sum_{\substack{i_{1},\ldots,i_{n},j_{1},\ldots,j_{n}=1
\\ $ \scriptsize{with} $ i_{1},\ldots,i_{n},j_{1},\ldots,j_{n}\\ $
\scriptsize{distinct} $
}}^{N}\delta_{(y_{i_{1}},\ldots,y_{i_{n}},y_{j_{1}},\ldots,y_{j_{n}})}.
\end{displaymath}
We can write
\begin{align*}
&\mathbb{E}\left(\langle f^{M},\tilde{\mu}_{t}^{N}\rangle^{2}\right)\\
&=\mathbb{E}\left(\int_{(0,\infty)^{2n}}P_{t}^{N}
f_{i_{1}}^{M}\left((x_{1},\ldots,x_{n})\right)
f_{i_{2}}^{M}\left((x_{n+1},\ldots,x_{2n})\right)\mu_{0}^{2n,N}(dx_{1},\ldots,dx_{2n})\right)\\
&\qquad+\circ(\frac{1}{N})
\end{align*}
Hence, by mimicking Lemma 2.2 we can write,
\begin{align*}
&\mathbb{E}\left(\langle f^{M},\tilde{\mu}_{t}^{N}\rangle^{2}\right)\\
&=\mathbb{E}\left(\int_{(0,\infty)^{2n}}P_{t}^{N}
f_{i_{1}}^{M}\left((x_{1},\ldots,x_{n})\right)
f_{i_{2}}^{M}\left((x_{n+1},\ldots,x_{2n})\right)\mu_{0}^{N}(dx_{1}),\ldots,\mu_{0}^{N}(dx_{2n})\right)\\
&\qquad+\circ(\frac{1}{N})
\end{align*}
By Theorem 3.4,
\begin{align*}
&P_{t}^{N}f_{i_{1}}^{M}\left((x_{1},\ldots,x_{n})\right)f_{i_{2}}^{M}\left((x_{n+1},\ldots,x_{2n})\right)\\
&\qquad\to
P_{t}f_{i_{1}}^{M}\left((x_{1},\ldots,x_{n})\right)f_{i_{2}}^{M}\left((x_{n+1},\ldots,x_{2n})\right)
\end{align*}
as $N\to\infty$ and is bounded by
$\|f_{i_{1}}^{M}\|_{\infty}\|f_{i_{2}}^{M}\|_{\infty}C^{2}t^{2n-2}$
where $C$ has been defined in 3.1. Hence by the Bounded Convergence
Theorem,

\begin{align*}\mathbb{E}\left(\langle f^{M},\tilde{\mu}_{t}^{N}
\rangle^{2}\right)\to
\int_{(0,\infty)^{2n}}&P_{t}f_{i_{1}}^{M}\left((x_{1},\ldots,x_{n})\right)
f_{i_{2}}^{M}\left((x_{n+1},\ldots,x_{2n})\right)\\
&\qquad\mu_{0}(dx_{1})\ldots\mu_{0}(dx_{2n})
\end{align*}

as $N\to\infty$. Now, by Theorem 3.4,
\begin{align*}&P_{t}f_{i_{1}}^{M}\left((x_{1},\ldots,x_{n})\right)f_{i_{2}}^{M}\left(
(x_{n+1},\ldots,x_{2n})\right)\\
&\qquad=P_{t}f_{i_{1}}^{M}\left((x_{1},\ldots,x_{n})\right)
P_{t}f_{i_{2}}^{M}\left((x_{n+1},\ldots,x_{2n})\right).
\end{align*}
Hence, as $N \to \infty,$
\begin{align*}
\mathbb{E}\left(\langle f^{M},\tilde{\mu}_{t}^{N}\rangle ^{2}\right)
\to
&\left(\int_{(0,\infty)^{n}}P_{t}f_{i_{1}}^{M}\left((x_{1},\ldots,x_{n})\right)
\mu_{0}(dx_{1})\ldots \mu_{0}(dx_{n})\right)\\
&\left(\int_{(0,\infty)^{n}}P_{t}f_{i_{2}}^{M}\left((x_{1},\ldots,x_{n})\right)
\mu_{0}(dx_{1})\ldots \mu_{0}(dx_{n})\right)\\
&=\langle f^{M},\tilde{\mu}_{t}\rangle ^{2}.
\end{align*}
\paragraph{Step 3:} We prove the Proposition 3.3.
We have proved before that as $M\to\infty,$ $\langle
f^{M},\tilde{\mu}_{t}\rangle\to\langle f,\tilde{\mu}_{t}\rangle$.
Hence,
\[\langle f^{M},\tilde{\mu}_{t}\rangle^{2}\to\langle f,\tilde{\mu}_{t}\rangle^{2}\] as $M \to \infty.$ Also, similarly to the way that we proved that
\[\lim_{M}\limsup_{N} \mathbb{E}\left(\left|\langle f(1-\Psi_{M}(m)),\tilde{\mu}_{t}^{N} \rangle\right|\right)=\lim_{M}\liminf_{N} \mathbb{E}\left(\left|\langle f(1-\Psi_{M}(m)),\tilde{\mu}_{t}^{N} \rangle\right|\right)=0\] we prove that
\[\lim_{M}\limsup_{N} \mathbb{E}\left(\left|\langle f(1-\Psi_{M}(m)),\tilde{\mu}_{t}^{N} \rangle\right|^{2}\right)=\lim_{M}\liminf_{N} \mathbb{E}\left(\left|\langle f(1-\Psi_{M}(m)),\tilde{\mu}_{t}^{N} \rangle\right|^{2}\right)=0.\]
Moreover, \begin{align*}&\mathbb{E}\left(\left|\langle
f\Psi_{M}(m),\tilde{\mu}_{t}^{N} \rangle \langle
f(1-\Psi_{M}(m)),\tilde{\mu}_{t}^{N} \rangle
\right|\right)^{2}\\
&\qquad\leq \mathbb{E}\left(\langle
f\Psi_{M}(m),\tilde{\mu}_{t}^{N}\rangle^{2}\right)\mathbb{E}\left(\langle
f(1-\Psi_{M}(m)),\tilde{\mu}_{t}^{N}\rangle^{2}\right). \end{align*}
Hence, we deduce that
\begin{align*}
&\lim_{M}\limsup_{N}\mathbb{E}\left(\left|\langle f\Psi_{M}(m),\tilde{\mu}_{t}^{N} \rangle \langle f(1-\Psi_{M}(m)),\tilde{\mu}_{t}^{N} \rangle \right|\right)^{2}\\
&\qquad=\lim_{M}\liminf_{N}\mathbb{E}\left(\left|\langle
f\Psi_{M}(m),\tilde{\mu}_{t}^{N} \rangle \langle
f(1-\Psi_{M}(m)),\tilde{\mu}_{t}^{N} \rangle \right|\right)^{2}=0.
\end{align*}
Hence, we obtain
\[\limsup_{N}\mathbb{E}\left(\langle f,\tilde{\mu}_{t}^{N} \rangle^{2}\right)=\liminf_{N}\mathbb{E}\left(\langle f,\tilde{\mu}_{t}^{N} \rangle^{2}\right)=\langle f,\tilde{\mu}_{t}\rangle^{2}.\]
Thus,
\[\lim_{N}\mathbb{E}\left(\langle f,\tilde{\mu}_{t}^{N} \rangle^{2}\right)=\langle f,\tilde{\mu}_{t}\rangle^{2}\] as required.
\subsection{Conclusion}
 For all $\tau \in \mathbb{T}$ and for all $f \in C_{b}\left(A_{\tau}(0,t)\right)$ we have proved that
\[\mathbb{E}\left(\langle f,\tilde{\mu}_{t}^{N}\rangle\right)\to \langle f,\tilde{\mu}_{t} \rangle \] and
\[\mathbb{E}\left(\langle f,\tilde{\mu}_{t}^{N}\rangle^{2}\right)\to \langle f,\tilde{\mu}_{t} \rangle ^{2}\]
as $N\to \infty$. So we deduce by the remark we did at the beginning
of the paper, that for all $f \in C_{b}\left(A_{\tau}(0,t)\right)$,
\[\langle f,\tilde{\mu}_{t}^{N} \rangle \to \langle f,\tilde{\mu}_{t} \rangle\]
as $N\to \infty$ in probability. Now let us prove our main result.

\section{Proof of Theorem 1.1}
Our aim in this section is to prove our main result stated in
Theorem 1.1, that is
\[\tilde{\mu}_{t}^{N}\to \tilde{\mu}_{t}\]
as $N \to \infty$ weakly in probability. In Section 3, we have
proved that for $\tau \in \mathbb{T}$ and for all $f \in
C_{b}\left(A_{\tau}(0,t)\right)$,
\[\langle f,\tilde{\mu}_{t}^{N} \rangle \to \langle f,\tilde{\mu}_{t} \rangle\]
as $N\to \infty$ in probability. To prove Theorem 1.1, we need to
prove that for all $f \in C_{b}\left(A(0,t)\right)$,
\[\langle f,\tilde{\mu}_{t}^{N} \rangle \to \langle f,\tilde{\mu}_{t} \rangle\]
as $N\to \infty$ in probability.
\subsection{A tightness argument}
We are reviewing here a particular case of the work of [1], [2].
The usual Marcus-Lushnikov process gives at each time the distribution in masses of the particles present in the system but does not retain any notion of configuration of these particles. We are going to define a process $(\overline{X}_{t}^{N})_{t\geq 0}$ (the Marcus Lushnikov on trees) on $\mathbb{T}(0,\infty)$, the space of trees on $(0,\infty)$. The space $\mathbb{T}(0,\infty)$ is given by
\begin{displaymath}
\mathbb{T}(0,\infty)=\bigcup_{\tau \in \mathbb{T}}\mathbb{T}_{\tau}(0,\infty)
\end{displaymath}
where $\mathbb{T}_{1}(0,\infty)=(0,\infty)$ and for $\tau=\{\tau_{1},\tau_{2}\}\in
\mathbb{T}$,
\[\mathbb{T}_{\tau}(0,\infty)=\big\{y=\{y_{1},y_{2}\}:y_{1}\in \mathbb{T}_{\tau_{1}}(0,\infty),
y_{2}\in \mathbb{T}_{\tau_{2}}(0,\infty)\big\}.\]
On $\mathbb{T}(0,\infty)$ define the mass function $m^{\prime}:\mathbb{T}(0,\infty) \to (0,\infty)$.
For $y \in \mathbb{T}_{1}(0,\infty)=(0,\infty)$, we set \[m^{\prime}(y)=y.\] Recursively for
$\tau=\{\tau_{1},\tau_{2}\} \in \mathbb{T}$, for
$y=\{y_{1},y_{2}\}\in \mathbb{T}_{\tau}(0,\infty),$ we set
\[m^{\prime}(y)=m^{\prime}(y_{1})+m^{\prime}(y_{2}).\]
Define also the counting function $n^{\prime}:\mathbb{T}(0,\infty) \to \mathbb{N}$ as follows.
For $y \in \mathbb{T}_{1}(0,\infty)=(0,\infty)$ we set
\[n^{\prime}(y)=1.\] Recursively for
$\tau=\{\tau_{1},\tau_{2}\} \in \mathbb{T}$, for
$y=\{y_{1},y_{2}\}\in \mathbb{T}_{\tau}(0,\infty),$ we set
\[n^{\prime}(y)=n^{\prime}(y_{1})+n^{\prime}(y_{2}).\]

Let us define $(\overline{X}_
{t})_{t\geq 0}$. Let $[N]=\{1,\ldots,N\}$ and let
$y_{1},\ldots,y_{N}>0$ be the masses associated to each particle in
$[N]$. Set
\begin{displaymath}
\overline{X}_{0}=\sum_{i=1}^{N}\delta_{y_{i}}.
\end{displaymath}
For each $i<j \in [N]$ take an independent random variable $T_{ij}$ such
that $T_{ij}$ is exponential with parameter $K\left(m^{\prime}(y_{i}),m^{\prime}(y_{j})\right)$, and
define
\begin{displaymath}
T=\min_{i<j}T_{ij}.
\end{displaymath}
Set $\overline{X}_{t}=\overline{X}_{0}$ for $t<T$ and
\[\overline{X}_{T}=\overline{X}_{t}-(\delta_{y_{i}}+\delta_{y_{j}}-\delta_{\{y_{i},y_{j}\}})\] if $T=T_{ij}$, then begin
the construction afresh from $\overline{X}_{T}$.
Let $(\overline{X}_{t}^{N})_{t\geq 0}$ be
Marcus-Lushnikov on $\mathbb{T}(0,\infty)$ with kernel $\frac{K}{N}$ starting from
\[
\overline{X}_{0}^{N}=\sum_{i=1}^{N}\delta_{y_{i}}.
\]
Set
\[\overline{\mu}_{t}^{N}=N^{-1}\overline{X}_{t}^{N}.\]
Let \[y_{t}:A(0,t)\to \mathbb{T}(0,\infty)\] be the map on forgetting times. Hence $\tilde{\mu}_{t}^{N}$ and $\overline{\mu}_{t}^{N}$ are related through the following equality
\begin{equation}
\overline{\mu}_{t}^{N}=\tilde{\mu}_{t}^{N}\circ y_{t}^{-1}.
\end{equation}
Define \[\tilde{\varphi}:\mathbb{T}(0,\infty)\to (0,\infty) \mbox{ by } \tilde{\varphi}=\varphi \circ m.\]
On $\mathbb{T}_{1}(0,\infty)=(0,\infty),$ $\tilde{\varphi}=\varphi$ and so
\[\langle \tilde{\varphi}^{2},\mu_{0}\rangle=\langle \varphi^{2},\mu_{0}\rangle <\infty\]
by (4). Then [1] tells us that we can find $n_{0}>0$ such that for all $N$,
\begin{displaymath}
\mathbb{P}\left(\overline{\mu}_{t}^{N}\left(\{y \in
\mathbb{T}(0,\infty):n^{\prime}(y)\geq
n_{0}\}\right)>C(\epsilon)\}\right)<\epsilon
\end{displaymath}
where $C(\epsilon)>0$.
Now by (22)
\[\overline{\mu}_{t}^{N}\left(\{y \in
\mathbb{T}(0,\infty):n^{\prime}(y)\geq
n_{0}\}\right)=\tilde{\mu}_{t}^{N}\left(\{\xi \in
A(0,t):n^{\prime}\left(y_{t}(\xi)\right)\geq
n_{0}\}\right).\]
Hence,
\begin{equation}\mathbb{P}\left(\tilde{\mu}_{t}^{N}\left(\{\xi \in
A(0,t):n\left(y_{t}(\xi)\right)\geq
n_{0}\}\right)>\frac{\epsilon}{\|f\|_{\infty}}\right)<\epsilon.
\end{equation}

\subsection{ Proof of Theorem 1.1} Take $f \in C_{b}\left(A(0,t)\right)$.
For $n \in \mathbb{N}$ consider
\begin{displaymath}
A_{n}=\bigcup_{\substack{\tau \in \mathbb{T}\\$ \scriptsize{with} $
n(\tau)\leq n}}A_{\tau}(0,t).
\end{displaymath}
Since \[A(0,t)=\bigcup_{\tau\in \mathbb{T}} A_{\tau}(0,t)\] it is
clear that
\begin{displaymath}
A(0,t)=\bigcup_{n\geq 1} A_{n}.
\end{displaymath}
Moreover, $\left(A_{n}\right)_{n \geq 1}$ forms an increasing
sequence. Set $f_{n}=f1_{A_{n}}$. Then, $f_{n}\in
C_{b}\left(A_{n}\right)$ and $f_{n}\to f$ as $n \to \infty$. We can
write
\begin{displaymath}
f_{n}=\sum_{\substack{\tau \in \mathbb{T}\\$ \scriptsize{with} $
n(\tau)\leq n}}f_{n}1_{A_{\tau}(0,t)}.
\end{displaymath}
Thus,
\begin{displaymath}
\langle f_{n},\tilde{\mu}_{t}^{N} \rangle=\sum_{\substack{\tau \in
\mathbb{T}\\$ \scriptsize{with} $ n(\tau)\leq n}} \langle
f_{n}1_{A_{\tau}(0,t)},\tilde{\mu}_{t}^{N} \rangle.
\end{displaymath}
Now, since $f_{n}1_{A_{\tau}(0,t)}\in
C_{b}\left(A_{\tau}(0,t)\right)$, by subsection 3.3,
\begin{displaymath}
\langle f_{n}1_{A_{\tau}(0,t)},\tilde{\mu}_{t}^{N} \rangle \to
\langle f_{n}1_{A_{\tau}(0,t)},\tilde{\mu}_{t} \rangle
\end{displaymath}
as $N\to\infty$, in probability. Thus, since the sum over
$\{\tau\in\mathbb{T} \mbox{ with } n(\tau)\leq n \}$ is finite, we
obtain
\begin{equation}
\langle f_{n},\tilde{\mu}_{t}^{N} \rangle \to \langle
f_{n},\tilde{\mu}_{t} \rangle
\end{equation}
as $N\to\infty$, in probability. Now let us prove our result. Also,
for all $n$, $f_{n}\to f$ as $n \to \infty$ and $|f_{n}|\leq |f|.$
Hence, by the dominated convergence theorem,
\begin{equation}
\tilde{\mu}_{t}(f_{n})\to \tilde{\mu}_{t}(f)
\end{equation}
as $n\to \infty$.

Now let us fix $\epsilon>0$. Consider,
\begin{align*}
|\langle f,\tilde{\mu}_{t}^{N} \rangle-\langle f,\tilde{\mu}_{t}
\rangle|\leq &|\langle f,\tilde{\mu}_{t}^{N} \rangle-\langle
f_{n},\tilde{\mu}_{t}^{N} \rangle|
+|\langle f_{n},\tilde{\mu}_{t}^{N} \rangle-\langle f_{n},\tilde{\mu}_{t} \rangle|\\
&\qquad+|\langle f_{n},\tilde{\mu}_{t} \rangle-\langle
f,\tilde{\mu}_{t} \rangle|
\end{align*}
By subsection 4.1, we know that we can find $n_{0}>0$, such that for all $N$,
\begin{equation}
\mathbb{P}\left(\tilde{\mu}_{t}^{N}\left(\{\xi \in
A(0,t):n^{\prime}\left(y_{t}(\xi)\right)\geq
n_{0}\}\right)>\frac{\epsilon}{\|f\|_{\infty}}\right)<\frac{\epsilon}{3}.
\end{equation}
Now, $f-f_{n}=f1_{A_{n}^{c}}$. Thus,
\[
|\langle f,\tilde{\mu}_{t}^{N} \rangle-\langle
f_{n},\tilde{\mu}_{t}^{N} \rangle|=|\langle
f-f_{n},\tilde{\mu}_{t}^{N} \rangle|\leq \|f\|_{\infty}|\langle
1_{A_{n}^{c}},\tilde{\mu}_{t}^{N} \rangle|
\]
For all $n \geq n_{0}, A_{n}^{c}\subseteq A_{n_{0}}^{c}$. So using
the relation (25) we obtain that
\begin{align*}
\mathbb{P}\left(|\langle f,\tilde{\mu}_{t}^{N} \rangle-\langle
f_{n},\tilde{\mu}_{t}^{N} \rangle|>\epsilon\right)&<
\mathbb{P}\left(\|f\|_{\infty}|\langle 1_{A_{n}^{c}},\tilde{\mu}_{t}^{N} \rangle|>\epsilon\right)\\
&<\frac{\epsilon}{3}
\end{align*}
for all $n \geq n_{0}$. Hence, for $n \geq n_{0}$, for all $N$,
\[\mathbb{P}\left(|\langle f,\tilde{\mu}_{t}^{N} \rangle-\langle f_{n},\tilde{\mu}_{t}^{N} \rangle|>\epsilon\right)<
\frac{\epsilon}{3}\]

For the third term, by 4.1, we can choose $n_{1}$ so that, for all
$n \geq n_{1}$,
\[|\langle f_{n},\tilde{\mu}_{t} \rangle-\langle f,\tilde{\mu}_{t} \rangle|<\frac{\epsilon}{3}\]

For the second term : Take $n \geq \max(n_{0},n_{1})$. By the
relation (24), we can find $N_{0}$ so that for all $N \geq N_{0}$,
\[
\mathbb{P}\left(|\langle f_{n},\tilde{\mu}_{t}^{N} \rangle-\langle
f_{n},\tilde{\mu}_{t} \rangle|>\epsilon\right)<\frac{\epsilon}{3}
\]
Hence, for $N>N_{0}$,
\[
\mathbb{P}\left(|\langle f,\tilde{\mu}_{t}^{N} \rangle-\langle
f,\tilde{\mu}_{t} \rangle|>\epsilon\right)<\epsilon
\]
that is
\[\langle f,\tilde{\mu}_{t}^{N} \rangle \to\langle f,\tilde{\mu}_{t} \rangle\] as $N\to \infty$ in probability.
Hence,
\[\tilde{\mu}_{t}^{N}\to \tilde{\mu}_{t}\]  as $N\to\infty$ weakly in probability as required.

\section{Appendix}
\subsection{A topology on $A(0,t)$}
We equip $A(0,t)$ with a topology. Define \[\hat{}:A(0,t)\to
A(0,1)\] as follows. For $\xi \in A_{1}(0,t)$ we set
$\hat{\xi}=\xi$. Recursively, for $\tau=\{\tau_{1},\tau_{2}\}\in
\mathbb{T},\xi=(s,\{\xi_{1},\xi_{2}\}) \in A_{\tau}(0,t)$ with
$\xi_{1} \in A_{\tau_{1}}(0,s)$ and $\xi_{2}\in A_{\tau_{2}}(0,s)$
we define $\hat{\xi}=(\frac{s}{t},\{\hat{\xi_{1}},\hat{\xi_{2}}\})
\in A_{\tau}(0,1)$ with $\hat{\xi_{1}} \in A_{\tau_{1}}(0,1)$ and
$\hat{\xi_{2}}\in A_{\tau_{2}}(0,1)$.

We are going to construct a topology on $A(0,1)$ and then rescaling we will obtain a topology on
$A(0,t)$. The set $\mathbb{T}$ is countable so we can give it a strict total order $<$. For $\tau=1\in \mathbb{T}, A_{1}(0,1)=(0,\infty)$. We
equip $ A_{1}(0,1)$ with the usual topology on $(0,\infty)$. For $\tau=\{\tau_{1},\tau_{2}\} \in \mathbb{T}$ with $\tau_{1}<\tau_{2}$, we can
identify $A_{\tau}(0,1)$ to be
\[A_{\tau}(0,1)=(0,1)\times A_{\tau_{1}}(0,1) \times A_{\tau_{2}}(0,1)\]

Now, for $\tau=\{1\}$,
\[A_{\{1\}}(0,1)=(0,1)\times A_{1}(0,1) \times A_{1}(0,1)\]
We equip $(0,1)$ with the usual topology on $\mathbb{R}$. We have
already given $A_{1}(0,1)$ a topology. Hence, we equip
$A_{\{1\}}(0,1)$ with the product topology (Tychonoff topology). By
induction, assume we have topologies on $A_{\tau_{1}}(0,1)$ and
$A_{\tau_{2}}(0,1)$ for $\tau_{1},\tau_{2}\in \mathbb{T}$.
Recursively for $\tau=\{\tau_{1},\tau_{2}\}\in \mathbb{T}$ with $\tau_{1}<\tau_{2}$, we equip
\[A_{\tau}(0,1)=(0,1)\times A_{\tau_{1}}(0,1) \times A_{\tau_{2}}(0,1)\]
with the product topology.
 Then ,
\[A(0,1)=\bigcup_{\tau\in \mathbb{T}}A_{\tau}(0,1)\] is naturally equipped with a topology. Finally, rescaling, we obtain a topology on $A(0,t)$ for $t \geq 0$. Similarly, we can equip $\tilde{A}(0,t)$ with a topology.

\subsection{A measure on $A_{i}^{y}(0,t)$}
Fix $0\leq t<T$. The set $\mathbb{T}_{n}^{\star}[n]$ is countable so we can give it a total order $<.$ Take
$i=\{i_{1},i_{2}\} \in \mathbb{T}_{n}^{\star}[n]$ with $i_{1}<i_{2}$. Let $y,y^{1}$ and $y^{2}$ by the vector of masses respectively associated to $i,i_{1}$ and
$i_{2}$. Our aim is to define a measure on $A_{i}^{y}(0,t)$. We are first going to construct a measure on $A_{i}^{y}(0,1)$ and then by a similar rescalling to the one done in subsection 5.1 we will obtain a measure on $A_{i}^{y}(0,t)$. Equip $A_{i}(0,1)$ with its Borel-$\sigma$ algebra $\mathcal{B}_{i}.$ We can
We can
identify $A_{i}^{y}(0,1)$ to be
\[A_{i}^{y}(0,1)=(0,1)\times A_{i_{1}}^{y^{1}}(0,1)\times A_{i_{2}}^{y^{2}}(0,1)\] as $i_{1}<i_{2}.$
For
$i=1,$ for $\xi \in
A_{i}^{y}(0,1),$ set
\[\nu_{i}^{y}(d\hat{\xi})=\delta_{y}.\]
Recursively for $i=\{i_{1},i_{2}\}
\in \mathbb{T}_{n}^{\star}[n]$ with $i_{1}<i_{2}$, with associated vector of masses $y,y^{1}$ and $y^{2}$, for
$\hat{\xi}=(s,\hat{\xi_{1}},\hat{\xi_{2}})\in
A_{i}(0,1)$, define
\[\nu_{i}^{y}(d\hat{\xi})=\nu_{i_{1}}^{y^{1}}(d\hat{\xi_{1}})\nu_{i_{2}}^{y^{2}}(d\hat{\xi_{2}})ds.\]
This defines a measure on the product space
\[\mathcal{B}\left((0,1)\right)\otimes
\mathcal{B}_{i_{1}}\otimes
\mathcal{B}_{i_{2}}\] where
$\mathcal{B}_{i_{1}}$ and
$\mathcal{B}_{i_{2}}$ are the respective Borel-$\sigma$
algebra on $A_{i_{1}}^{y^{1}}(0,1)$ and
$A_{i_{2}}^{y^{2}}(0,1)$. Then rescalling as in 5.1 we obtain a
measure on $A_{i}^{y}(0,t)$.

\subsection{A limit measure on $\tilde{A}(0,t)$}
Let $E=(0,\infty).$ Fix $0\leq t<T$. Fix $n>0$. We define a limit measure on $\tilde{A}(0,t)$ as follows.
For $\xi \in A_{k}(0,t)$ with $k\in\mathbb{N}$, set
\[\tilde{\mu}_{t}^{\prime}(d\xi)=\exp\left(-\int_{0}^{t}\int_{E}K(y,y^{\prime})\mu_{r}(dy^{\prime})\mathrm{d}r\right)\mu_{0}(dy)\]
where $y=m(\xi)$ and $(\mu_{r})_{r<T}$ is the strong
deterministic solution to the generalized Smoluchowski equation (3).
For $i=\{i_{1},i_{2}\}\in
\mathbb{T}_{n}^{\star}\mathbb{N}$ with type $\tau(i)=\tau=\{\tau_{1},\tau_{2}\} \in \mathbb{T}$ , for $\xi=(s,\{\xi_{1},\xi_{2}\})\in
A_{i}(0,t)$, with $s<t$, define recursively
\begin{equation} \label{eq2345}
\tilde{\mu}_{t}^{\prime}(d\xi)=\epsilon(\tau)K(m(\xi_{1}),m(\xi_{2}))\tilde{\mu}_{s}^{\prime}(d\xi_{1})\tilde{\mu}_{s}^{\prime}(d\xi_{2})\exp\left(-\int_{s}^{t}\int_{E}K(y,y^{\prime})\mu_{r}(dy^{\prime})\mathrm{d}r\right)ds
\end{equation} where $y=m(\xi)$ and $\epsilon(\tau)=1$ if $\tau_{1}\neq\tau_{2}$ and $\epsilon(\tau)=\frac{1}{2}$ if $\tau_{1}=\tau_{2}.$

For $\xi \in A_{k}(0,t)$ with $k\in\mathbb{N}$, set
\[K_{\xi}=1\]

For $i=\{i_{1},i_{2}\}\in
\mathbb{T}_{n}^{\star}\mathbb{N}$, $\xi=(s,\{\xi_{1},\xi_{2}\})\in
A_{i}(0,t)$, define recursively
\[K_{\xi}=K(m(\xi_{1}),m(\xi_{2}))K_{\xi_{1}}K_{\xi_{2}}\]

\begin{theorem} Let $i=\{i_{1},i_{2}\}\in \mathbb{T}_{n}^{\star}\mathbb{N}$ with type $\tau=\tau(i) \in \mathbb{T}$.
\begin{enumerate}
\item Then for $\xi=(s,\{\xi_{1},\xi_{2}\})\in A_{i}(0,t)$,
\[\tilde{\mu}_{t}^{\prime}(d\xi)=2^{-q(\tau)}K_{\xi}\exp\left(-\int_{\Delta(\xi)}\int_{E}K(y_{r},y^{\prime})\mu_{\Pi(r)}(dy^{\prime})dr\right)\nu_{i}^{y}(d\xi) \mu_{0}(dy_{1})\ldots\mu_{0}(dy_{n})\] where $y=(y_{1},\ldots,y_{n})$ is the vector of masses associated to $\xi$ and $\nu_{i}^{y}$ is the measure described in subsection 5.2.
\item Take $f\in C_{b}(A_{\tau}(0,t)$. Let $f_{i}=f\circ g_{i}$ where $g_{i}:A_{i}(0,t)\to A_{\tau}(0,t)$ is the map on forgetting labels. Then,\[\int_{A_{\tau}(0,t)}f(\xi)\tilde{\mu}_{t}(d\xi)=\int_{A_{i}(0,t)}f_{i}(\xi)\int_{A_{\tau}(0,t)}f(\xi)\tilde{\mu}_{t}(d\xi)\]
\end{enumerate}
\end{theorem}

{\it Proof of Theorem 5.1}: ~~\\~~\\
(a) Let us do it by induction. It is clearly
true for $i\in \mathbb{T}_{1}^{\star}\mathbb{N}$. Fix $n>0.$ Suppose
it is true for all $i \in
\mathbb{T}_{k}^{\star}\mathbb{N}$ with $k\leq n-1$. Is is true for
$k= n$? Take $i=\{i_{1},i_{2}\}\in
\mathbb{T}_{n}^{\star}\mathbb{N}$. Then, writing
$n(i_{1})=k$ and $n(i_{2})=n-k$ we have
\[n(i_{1})\leq n-1\] and \[n(i_{2})\leq n-1\]
Take $\xi=(s,\{\xi_{1},\xi_{2}\})\in A_{i}(0,t)$. Without
loss of generality, assume that $i_{1}$ is formed from
the particle $1,\ldots,k$ with associated masses
$y^{1}=(y_{1},\ldots,y_{k})$ and that $i_{2}$ is formed
from the particle $k+1,\ldots,n$ with associated masses
$y^{2}=(y_{k+1},\ldots,y_{n})$. Using the induction hypothesis, we
have,
\begin{align*}
\tilde{\mu}_{s}^{\prime}(d\xi_{1})=&2^{-q(\tau_{1})}K_{\xi_{1}}\exp\left(-\int_{\Delta(\xi_{1})}\int_{E}K(y_{r},y^{\prime})
\mu_{\Pi(r)}(dy^{\prime})dr\right)\\&\nu_{i_{1}}^{y^{1}}(d\xi_{1})
\mu_{0}(dy_{1})\ldots\mu_{0}(dy_{k})
\end{align*}
and
\begin{align*}
\tilde{\mu}_{s}^{\prime}(d\xi_{2})=&2^{-q(\tau_{2})}K_{\xi_{2}}\exp\left(-\int_{\Delta(\xi_{2})}\int_{E}K(y_{r},y^{\prime})
\mu_{\Pi(r)}(dy^{\prime})dr\right)\\&\nu_{i_{2}}^{y^{2}}(d\xi_{2})
\mu_{0}(dy_{k+1})\ldots\mu_{0}(dy_{n})
\end{align*}
Using the relation $(27)$, we have
\begin{align*}
\tilde{\mu}_{t}^{\prime}(d\xi)=&\epsilon(\tau)2^{-q(\tau_{1})}2^{-q(\tau_{2})}K(m(\xi_{1}),m(\xi_{2}))K_{\xi_{1}}K_{\xi_{2}} \nu_{i_{1}}^{y^{1}}(d\xi_{1})\nu_{i_{2}}^{y^{2}}(d\xi_{2})ds\mu_{0}(dy_{1})\ldots\mu_{0}(dy_{n})\\
&\exp\left(-\int_{\Delta(\xi_{1})}\int_{E}K(y_{r},y^{\prime})\mu_{\Pi(r)}(dy^{\prime})dr\right)\exp\left(-\int_{\Delta(\xi_{2})}\int_{E}K(y_{r},y^{\prime})\mu_{\Pi(r)}(dy^{\prime})dr\right)\\&\exp\left(-\int_{s}^{t}\int_{E}K(y,y^{\prime})\mu_{r}(dy^{\prime})\mathrm{d}r\right)
\end{align*}
Now, $K_{\xi}=K(m(\xi_{1}),m(\xi_{2}))K_{\xi_{1}}K_{\xi_{2}}$,
$\nu_{i}^{y}(d\xi)=\nu_{i_{1}}^{y^{1}}(d\xi_{1})\nu_{i_{2}}^{y^{2}}(d\xi_{2})ds$,
$\Delta(\xi)=\Delta(\xi_{1})\cup\Delta(\xi_{2})\cup
\left((s,t]\times{\xi}\right)$ and $2^{-q(\tau)}=2^{-q(\tau_{1})}2^{-q(\tau_{2})}2^{-\epsilon(\tau)}$. Hence,
\[\tilde{\mu}_{t}^{\prime}(d\xi)=2^{-q(\tau)}K_{\xi}\exp\left(-\int_{\Delta(\xi)}\int_{E}K(y_{r},y^{\prime})\mu_{\Pi(r)}(dy^{\prime})dr\right)\nu_{i}^{y}(d\xi) \mu_{0}(dy_{1})\ldots\mu_{0}(dy_{n})\] as required.
\begin{flushright}
$\square$
\end{flushright}

\newpage
\subsection{Some simulations}

Here are some simulations under Visual Basic of the Marcus-Lushnikov process on trees.
The graphics below represent trees that have been simulated following the Marcus-Lushnikov process on trees with different kernel $K$ and an initial number of particles $N$. In these simulations, all the initial particles have mass $1$. These pictures show for each kernel the sort of trees limit we can expect to find in the limit measure.

\begin{figure}[htp]
\centering
\includegraphics[scale=0.45]{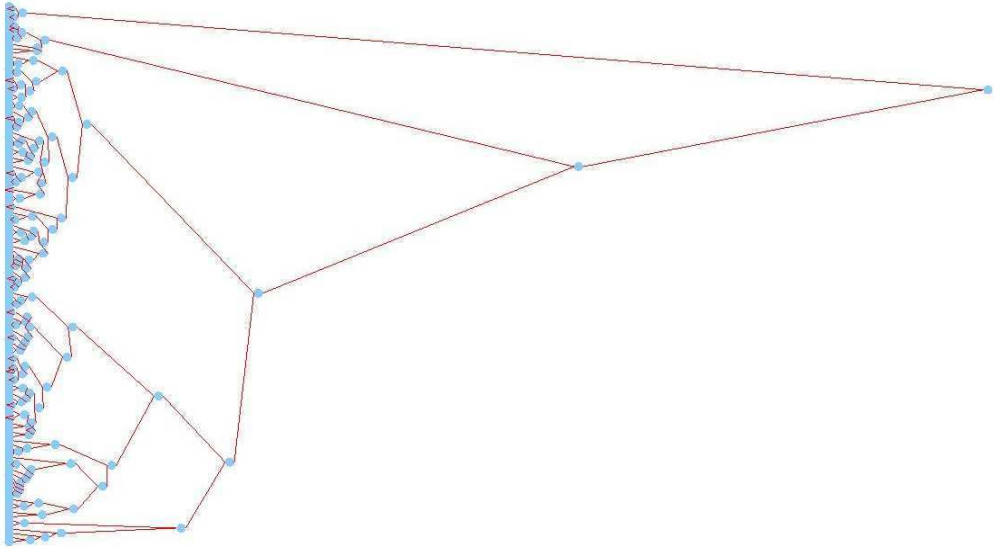}
\caption{$N=128$, $K(x,y)=1$}
\end{figure}

\begin{figure}[htp]
\centering
\includegraphics[scale=0.45]{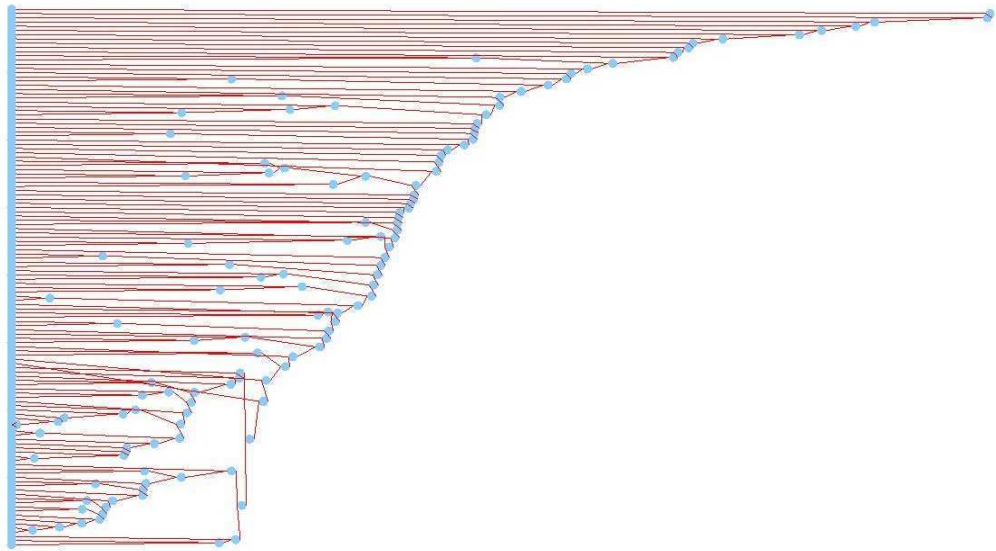} 
\caption{$N=128$, $K(x,y)=xy$}
\end{figure}

\begin{figure}[htp]
\centering
\includegraphics[scale=0.45]{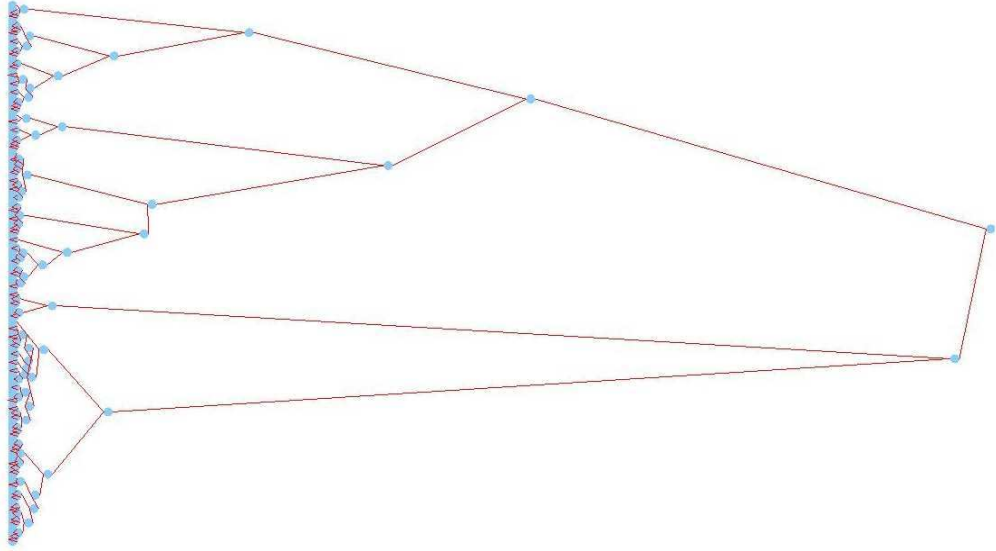}
\caption{$N=128$, $K(x,y)=1/(x+y+1)$}
\end{figure}

\end{document}